\documentclass[11pt]{amsart}
\usepackage[dvipsnames,usenames]{color}
\usepackage{hyperref}
\usepackage{graphicx}
\usepackage{epsfig}
\usepackage[latin1]{inputenc}
\usepackage{amsmath}
\usepackage{amsfonts}
\usepackage{amssymb}
\usepackage{amsthm}
\usepackage{amscd}
\usepackage{verbatim}
\usepackage{subfigure}
\usepackage{caption}
\usepackage{pinlabel}
\usepackage{stmaryrd}
\usepackage{enumerate, enumitem}
\usepackage{todonotes}
\usepackage{bm}
\usepackage{thmtools}
\usepackage{thm-restate}
\usepackage{lipsum}
\usepackage{setspace}
\usepackage{mathtools}
\usepackage[all]{xypic}
\usepackage[abs]{overpic}

\allowdisplaybreaks

\usepackage{mathdots}
\usepackage{tikz}
\usepackage{tikz-cd}
\usetikzlibrary{arrows}
\usetikzlibrary{decorations.pathreplacing}
\usepackage{verbatim}
\usetikzlibrary{cd}
\tikzset{taar/.style={double, double equal sign distance, -implies}}
\tikzset{amar/.style={->, dotted}}
\tikzset{dmar/.style={->, dashed}}
\tikzset{aar/.style={->, very thick}}

\usepackage{placeins}
\newtheorem{theorem}{Theorem}[section]

\theoremstyle{definition}
\newtheorem{definition}[theorem]{Definition}

\theoremstyle{remark}
\newtheorem{example}[theorem]{Example}

\def\F{\mathbb{F}}
\def\N{\mathbb{N}}

\def\Z{\mathbb{Z}}

\def\del{\partial}

\def\CFKi{\CFK^{\infty}}

\def\HFK {\mathit{HFK}}
\newcommand\HFKhat{\widehat{\HFK}}

\def\CFK{\mathit{CFK}}

\author[A. Antal]{Anna Antal}
\email{anna.antal@yale.edu}

\author[S. Pritchard]{Sarah Pritchard}
\email{sarah.pritchard@gatech.edu}

\numberwithin{equation}{section}

\title{A note on the involutive concordance invariants for certain (1,1)-knots}

\begin{document}
\maketitle
\begin{abstract} We compute the involutive concordance invariants for the 10- and 11-crossing (1,1)-knots. \end{abstract}
\color{black}
\tableofcontents
\section{Introduction}
Heegaard Floer homology is a suite of invariants of 3-manifolds, knots, and links introduced in the early 2000's by P. Ozsv\'{a}th and Z. Szab\'{o} \cite{OSknots,OSsurvey}, and in the knot case independently by J. Rasmussen \cite{RasmussenThesis}. In the knot variant, Heegaard Floer homology associates to a knot $K$ a $\Z \oplus \Z$-filtered, $\Z$-graded chain complex over $\F_2[U,U^{-1}]$ called $\CFKi{(K)}.$ Many classical knot invariants can be recovered from this chain complex. For example, $\CFKi{(K)}$ contains the data of the Alexander polynomial \cite{OSknots}, the knot genus \cite{OS:four}, and whether a knot is fibred \cite{Ghigginifibred,Nifibred}. 

We work with involutive Heegaard Floer homology, developed by K. Hendricks and C. Manolescu  in 2015 \cite{HM:involutive} as a refinement to Heegaard Floer homology. In the knot version, involutive Heegaard Floer homology additionally considers a skew-filtered automorphism $$\iota_K: \CFKi{(K)} \rightarrow \CFKi{(K)},$$ which is order four up to filtered chain homotopy. Using this supplementary information, involutive Heegaard Floer homology introduced two new knot concordance invariants, $\underline{V}_0(K)$ and $\overline{V}_0(K)$, which are variants of an existing concordance invariant $V_0(K)$ from Heegaard Floer homology. The involutive concordance invariants are interesting because unlike other concordance invariants arising from Heegaard Floer homology such as $\tau, \epsilon$, and $\nu$, they do not necessarily vanish on knots of finite concordance order such as the figure-eight knot \cite{HM:involutive}. Hendricks and Manolescu give combinatorial computations of the involutive concordance invariants for $L$-space knots (which include the torus knots) and thin knots (which include alternating and quasi-alternating knots).

Following a similar, but slightly more refined, strategy to Hendricks and Manolescu's computation for thin knots \cite[Section 8]{HM:involutive}, in this note we compute the involutive concordance invariants $\underline{V}_0(K)$ and $\overline{V}_0(K)$ for all 10 and 11-crossing (1,1)-knots that are neither $L$-space nor thin. We also include the values of the ordinary concordance invariant $V_0(K)$ for each knot, for easy comparison. Our computations appear in Section \ref{sec:results}. Taken together they confirm the following.

\begin{theorem}
For any 10- or 11-crossing (1,1)-knot $K$, the skew-filtered chain homotopy equivalence class of $\iota_K$, and therefore the value of the concordance invariants $\underline{V}_0(K)$ and $\overline{V}_0(K)$, is determined by $\CFKi{(K)}.$
\end{theorem}

\subsection*{Organization} This paper is organized as follows.  In Section \ref{sec:background} we briefly introduce Heegaard Floer homology and involutive Heegaard Floer homology; specifically, in Section \ref{sub:cfki} we review the definition of the chain complex $\CFKi{(K)}$ and the concordance invariant $V_0(K)$, and in Section \ref{sec:involutive} we summarize the properties of the map $\iota_K$ and the definitions of the involutive concordance invariants $\underline{V}_0(K)$ and $\overline{V}_0(K)$. In Section \ref{sec:1,1} we discuss parameterizations of (1,1)-knots and Heegaard diagrams. We provide some example computations of the involutive concordance invariants in Section \ref{sec:examples}; Section \ref{sub:11n57ex} outlines the computation for the knot $11n_{57}$, and Section \ref{sub:10_161ex} outlines the computation for the knot $10_{161}$. Then, in Section \ref{sec:results} we provide a list of the values of the involutive concordance invariants for the 10- and 11-crossing (1,1)-knots for which they were not previously known, along with some amendments to the literature.

\subsection*{Acknowledgments} This project was carried out during the Rutgers DIMACS REU in Summer 2021; we thank the organizers of the REU for their support and Kristen Hendricks and Karuna Sangam for supervising this project. We are also grateful to Zipei Nie for sending us his Python code to compute the knot Floer homology of $(1,1)$ knots. Finally, we thank the anonymous referee for helpful comments and corrections. Both authors were supported by NSF CAREER Grant DMS-2019396.

\section{Heegaard Floer homology and involutive Heegaard Floer homology} 
\label{sec:background}

\subsection{The chain complex $CFK^{\infty}(K)$}
\label{sub:cfki}

We now introduce the complex $CFK^{\infty}(K)$ abstractly; in Section \ref{sec:1,1} we will go over its construction in a special case. To a knot $K$, Heegaard Floer homology \cite{OSsurvey,OSknots} associates a ($\Z \oplus \Z$)-filtered, $\Z$-graded chain complex $CFK^{\infty}(K)$ over $\F_2[U,U^{-1}]$. (Strictly speaking, the construction of $CFK^{\infty}(K)$ involves some choices which produce chain complexes which are chain homotopy equivalent via canonical chain homotopies; here and throughout, we will take some model for this chain homotopy equivalence class.)

For the mirror image $\overline{K}$ of a knot $K$, $CFK^{\infty}(\overline{K}) = \textrm{Hom}_{\F_2[U,U^{-1}]}(CFK^{\infty}(K),\F_2[U,U^{-1}])$, the dual of $CFK^{\infty}(K)$ over the field $\F_2[U,U^{-1}].$ We can describe the horizontal and vertical components of $\partial$ in following way: 

\begin{definition}
Decompose the differential $\partial$ on $\CFK^{\infty}(K)$ as $$\partial = \sum_{i,j \in \N} \partial_{ij},$$ where $\partial_{ij}$ lowers the horizontal grading by $i$ and the vertical grading by $j$. Then the horizontal and vertical differentials are
$$\partial_{\mathrm{horz}} = \sum_{i} \partial_{i0}$$ and
$$\partial_{\mathrm{vert}} = \sum_{j} \partial_{0j}.$$ 
\end{definition}

\noindent Given $S \subseteq \Z \oplus \Z$, we let $C\{S\} \subseteq CFK^{\infty}(K)$ denote the set of elements with planar gradings in $S$; if this is closed under $\partial$ it is a subcomplex of $CFK^{\infty}(K)$. The following complexes are important examples:

\begin{definition}
The quotient complexes $C\{i=0\}$ and $C\{j=0\}$ are
\[C\{i=0\} = C\{(i,j): i\leq 0\} / C\{(i,j):i < 0\} \] and
\[C\{j=0\} = C\{(i,j): j\leq 0\} / C\{(i,j):j < 0\}. \]
\end{definition}

\noindent Since both $C\{(i,j): i\leq 0\}$ and $C\{(i,j):i < 0\}$ are subcomplexes of $CFK^{\infty}(K)$, their quotient is a chain complex. In the same way, $C\{j=0\}$ is a subcomplex also. The Euler characteristic of the associated graded of $C\{i=0\}$ for a knot $K$ is the Alexander polynomial $\Delta_K(t)$. Additionally, the homology of the subcomplex $C\{i=0\}$ equipped with the boundary map $\partial_{\textrm{vert}}$ is
\[H_*\big(C\{i=0\}, \partial_{\textrm{vert}} \big) = \F_{(0)}.\] Likewise,
\[H_*\big(C\{j=0\}, \partial_{\textrm{horz}} \big) = \F_{(0)}.\]

\noindent We are furthermore interested in the subcomplex $A_0^-$, which has the following definition:

\begin{definition}
The subcomplex $A_0^-$ is defined by 
$$A_0^- = C\{(i,j): i,j \leq 0\}.$$
\end{definition}

Now, we consider some conventions for the properties of elements in $\CFKi{(K)}$. If $x \in \CFKi{(K)}$ is in grading $(i,j)$, then $Ux$ is in grading $(i-1,j-1)$. There are two other important gradings associated to an element of $\CFKi{(K)}$.

\begin{definition}
Suppose that $\CFKi{(K)}$ is normalized such that the elements $x = U^0x$ lie in planar gradings $(0,j)$, i.e. on the $j$-axis. Then, the \textit{Alexander grading} $A(x)$ of an element $x \in \CFKi{(K)}$ is the $j$-grading of $x$.

The \textit{homological grading} $M(x)$ of an element $x \in \CFKi{(K)}$ is determined by the following conventions. The element of $\CFKi{(K)}$ that generates $H_*\big(C\{i=0\}, \partial_{\textrm{vert}}\big)$ has homological grading 0. The boundary map $\partial$ lowers the homological grading of an element by 1. Multiplication by U lowers the homological grading of an element by 2.
\end{definition}

\noindent There are many concordance invariants that can be derived from $CFK^{\infty}(K)$, including the invariant $V_0(K)$, which is defined as follows. 

\begin{definition}
The concordance invariant $V_0(K)$ for a knot $K$ is given by
$$V_0(K) = -\frac{1}{2}\textrm{ max }\{r: \exists x \in H_r(A_0^-) \textrm{ such that } U^nx \neq 0 \textrm{ for all n }\}.$$ Here $r$ is the homological grading.
\end{definition}

\begin{example}
\label{ex:rht}
The complex $CFK^{\infty}(K)$, where $K$ is the right-handed trefoil knot $3_1$, is shown in Figure \ref{fig:rht}. The associated graded of $C\{i=0\}$ has Euler characteristic
\[\chi (C\{i=0)\}) = \sum_{j \in \Z} \chi(C_{0,j})t^j = t^{-1} - 1 + t.\] So, the Alexander polynomial of the right handed trefoil is $\Delta_{3_1}(t) = t^{-1} - 1 + t$. The homology of $A_0^-$ is $\F_2[U] \langle [Ua] \rangle$. Thus, $V_0(3_1)$ is $-\frac{1}{2}$ times the homological grading of $Ua$, which yields $$V_0(3_1) = -\frac{1}{2}(-2) = 1.$$

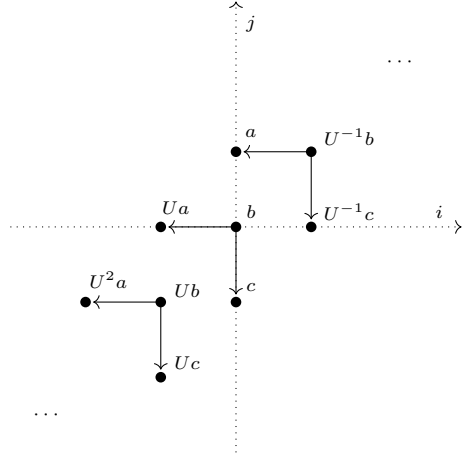
\begin{figure}
\begin{center}
    \begin{tikzpicture}\tikzstyle{every node}=[font=\tiny] 
    \path[->][dotted](0,-3)edge(0,3);
    \path[->][dotted](-3,0)edge(3,0);
    \node(1) at (2.2,2.2){$\cdots$};
    \node(2) at (-2.5,-2.5){$\cdots$};
    \node(3) at (2.7,.2){$i$};
    \node(4) at (.2,2.7){$j$};
    
    \node(5)at (0.2,1.25){$a$};
    \node(6) at (0.2,0.2){$b$};
    \node(7) at (.2,-.8){$c$};
    \node(8) at (-0.8,0.25){$Ua$};
    \node(9) at (-.65,-.85){$Ub$};
    \node(10) at (-.65,-1.8){$Uc$};
    \node(11) at (-1.7,-.7){$U^2a$};
    \node(12) at (1.5,1.2){$U^{-1}b$};
    \node(13) at (1.5,0.2){$U^{-1}c$};
    
    \fill(1,1)circle [radius=2pt];
    \fill(1,0)circle [radius = 2pt];
    \fill(0,1)circle [radius =2pt];
    \path[->](1,1) edge (0.1,1);
    \path[->](1,1) edge (1,0.1);
    
    \fill(0,0)circle [radius=2pt];
    \fill(0,-1)circle [radius = 2pt];
    \fill(-1,0)circle [radius =2pt];
    \path[->](0,0) edge (-0.9,0);
    \path[->](0,0) edge (0,-0.9);
    
    \fill(-1,-1)circle [radius=2pt];
    \fill(-1,-2)circle [radius = 2pt];
    \fill(-2,-1)circle [radius =2pt];
    \path[->](-1,-1) edge (-1.9,-1);
    \path[->](-1,-1) edge (-1,-1.9);
    
    \end{tikzpicture}
    \caption{The chain complex $\CFKi{(3_1)}$ has three generators, which we call $a,b,$ and $c$. The arrows in the picture represent the boundary map $\partial$. The homological grading of $a$ is 0; $b$ has homological grading -1; and $c$ has homological grading -2.}
    \label{fig:rht}
\end{center}
\end{figure}

\end{example}
\label{ex:lht}

\begin{example}
We now consider the left-handed trefoil knot, which is the mirror of the right-handed trefoil knot. $CFK^{\infty}(\overline{3_1})$ is shown in Figure \ref{fig:lht}. The associated graded of $C\{i=0\}$ has Euler characteristic
\[\chi (C\{i=0)\}) = \sum_{j \in \Z} \chi(C_{0,j})t^j = t^{-1} - 1 + t.\] So, the Alexander polynomial of the left handed trefoil is $\Delta_{\overline{3_1}}(t) = t^{-1} - 1 + t$. The homology of $A_0^-$ is $\F_2[U] \langle [Ua+c] \rangle$. Thus, $V_0(\overline{3_1})$ is $-\frac{1}{2}$ times the homological grading of $Ua+c$, which yields $$V_0(\overline{3_1}) = -\frac{1}{2}(0) = 0.$$

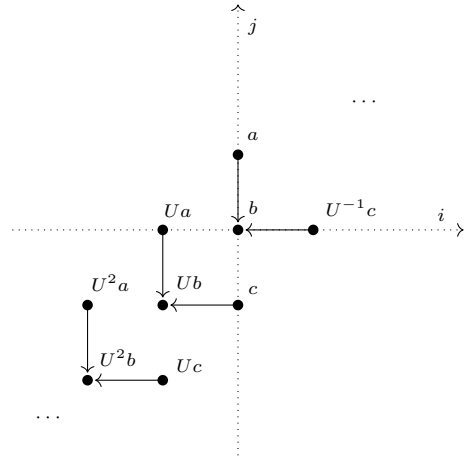
\begin{figure}
    \centering
    \begin{tikzpicture}\tikzstyle{every node}=[font=\tiny] 
    \path[->][dotted](0,-3)edge(0,3);
    \path[->][dotted](-3,0)edge(3,0);
    \node(1) at (1.7,1.7){$\cdots$};
    \node(2) at (-2.5,-2.5){$\cdots$};
    \node(3) at (2.7,.2){$i$};
    \node(4) at (.2,2.7){$j$};
    
    \node(5)at (0.2,1.25){$a$};
    \node(6) at (0.2,0.3){$b$};
    \node(7) at (.2,-.8){$c$};
    \node(8) at (-0.8,0.25){$Ua$};
    \node(9) at (-.65,-.7){$Ub$};
    \node(10) at (-.65,-1.8){$Uc$};
    \node(11) at (-1.7,-.7){$U^2a$};
    \node(12) at (-1.6,-1.7){$U^{2}b$};
    \node(13) at (1.5,0.3){$U^{-1}c$};
    
    \fill(1,0)circle [radius = 2pt];
    \fill(0,1)circle [radius =2pt];
    \path[->](1,0) edge (0.1,0);
    \path[->](0,1) edge (0,0.1);
    
    \fill(0,0)circle [radius=2pt];
    \fill(0,-1)circle [radius = 2pt];
    \fill(-1,0)circle [radius =2pt];
    \path[->](-1,0) edge (-1,-0.9);
    \path[->](0,-1) edge (-0.9,-1);
    
    \fill(-1,-1)circle [radius=2pt];
    \fill(-1,-2)circle [radius = 2pt];
    \fill(-2,-1)circle [radius =2pt];
    \fill(-2,-2)circle [radius =2pt];
    \path[->](-2,-1) edge (-2,-1.9);
    \path[->](-1,-2) edge (-1.9,-2);
    
    \end{tikzpicture}
    \caption{The chain complex $CFK^{\infty}(\overline{3_1})$ has three generators, which we call $a,b,$ and $c$. The homological grading of $a$ is 2; $b$ has homological grading 1; and $c$ has homological grading 0.}
    \label{fig:lht}
\end{figure}
\end{example}
\begin{example}
\label{ex:figeight}
Now we consider $CFK^{\infty}(K)$, where $K$ is the figure-eight knot $4_1$, as shown in Figure \ref{fig:figEight}. The associated graded of $C\{i=0\}$ has Euler characteristic
\[\chi (C\{i=0)\}) = \sum_{j \in \Z} \chi(C_{0,j})t^j = t^{-1} -3 + t.\] So, the Alexander polynomial of the figure eight knot is $\Delta_{4_1}(t) = t^{-1} -3 + t.$ The homology of $A_0^-$ is $\F_2[U] \langle [x] \rangle \cup [e]$. Thus, $V_0(4_1)$ is $-\frac{1}{2}$ times the homological grading of $x$, which yields $$V_0(4_1) = -\frac{1}{2}(0) = 0.$$

\begin{figure}
\begin{center}
    \begin{tikzpicture}\tikzstyle{every node}=[font=\tiny] 
    \path[->][dotted](0,-3.75)edge(0,4);
    \path[->][dotted](-4,0)edge(4,0);
    \fill(0,0)circle [radius=2pt];
    \fill(-.15,-.15)circle[radius=2pt];
    \fill(-.15,-1.15)circle [radius = 2pt];
    \fill(-1.15,-.15)circle [radius =2pt];
    \path[->](-.15,-.15)edge(-.15,-1.08);
    \path[->](-.15,-.15)edge(-1.08,-.15);
    \fill(-1.15,-1.15)circle[radius = 2pt];
    \path[->](-1.15,-.15)edge(-1.15,-1.08);
    \path[->](-.15,-1.15)edge(-1.08,-1.15);
    \begin{scope}[shift={(-1.3,-1.3)}]
    \fill(0,0)circle [radius=2pt];
    \fill(-.15,-.15)circle[radius=2pt];
    \fill(-.15,-1.15)circle [radius = 2pt];
    \fill(-1.15,-.15)circle [radius =2pt];
    \path[->](-.15,-.15)edge(-.15,-1.08);
    \path[->](-.15,-.15)edge(-1.08,-.15);
    \fill(-1.15,-1.15)circle[radius = 2pt];
    \path[->](-1.15,-.15)edge(-1.15,-1.08);
    \path[->](-.15,-1.15)edge(-1.08,-1.15);
    \node (1) at (.3,-.2){$Ux$};
    \node(2) at (-.42,-.38){$Ua$};
    \node(3) at (-1.5,-.1){$Ub$};
    \node(4) at (-.7,-.98){$U^2e$};
    \node(5) at (.2,-1.25){$Uc$};
    \end{scope}
    \begin{scope}[shift={(1.3,1.3)}]
    \fill(0,0)circle [radius=2pt];
    \fill(-.15,-.15)circle[radius=2pt];
    \fill(-.15,-1.15)circle [radius = 2pt];
    \fill(-1.15,-.15)circle [radius =2pt];
    \path[->](-.15,-.15)edge(-.15,-1.08);
    \path[->](-.15,-.15)edge(-1.08,-.15);
    \fill(-1.15,-1.15)circle[radius = 2pt];
    \path[->](-1.15,-.15)edge(-1.15,-1.08);
    \path[->](-.15,-1.15)edge(-1.08,-1.15);
    \end{scope}
    \begin{scope}[shift={(2.6,2.6)}]
    \fill(0,0)circle [radius=2pt];
    \fill(-.15,-.15)circle[radius=2pt];
    \fill(-.15,-1.15)circle [radius = 2pt];
    \fill(-1.15,-.15)circle [radius =2pt];
    \path[->](-.15,-.15)edge(-.15,-1.08);
    \path[->](-.15,-.15)edge(-1.08,-.15);
    \fill(-1.15,-1.15)circle[radius = 2pt];
    \path[->](-1.15,-.15)edge(-1.15,-1.08);
    \path[->](-.15,-1.15)edge(-1.08,-1.15);
    \end{scope}
    \fill(2.75,2.75)circle[radius=2pt];
    \fill(-2.6,-2.6)circle[radius=2pt];
    \fill(-2.75,-2.75)circle[radius=2pt];
    \node (1) at (.2,-.2){$x$};
    \node(2) at (-.35,-.35){$a$};
    \node(3) at (-1.35,-.2){$b$};
    \node(4) at (-.8,-.95){$Ue$};
    \node(5) at (-.2,-1.4){$c$};
    \node(6) at (.35,.3){$e$};
    \node(7) at (.48,1.4){$U^{-1}b$};
    \node(8) at (1.65,.25){$U^{-1}c$};
    \node(9) at (3.2,3.2){$\cdots$};
    \node(10) at (-3.2, -3.2){$\cdots$};
    \node(11) at (.2,3.7){$j$};
    \node(12) at (3.7,.2){$i$};
    \end{tikzpicture}
    \caption{The chain complex $CFK^{\infty}(4_1)$ has five generators, which we call $a,b,c,e$ and $x$. The homological gradings of $a,e$ and $x$ are 0; $b$ has homological grading 1; and $c$ has homological grading -1.}
    \label{fig:figEight}
\end{center}
\end{figure}
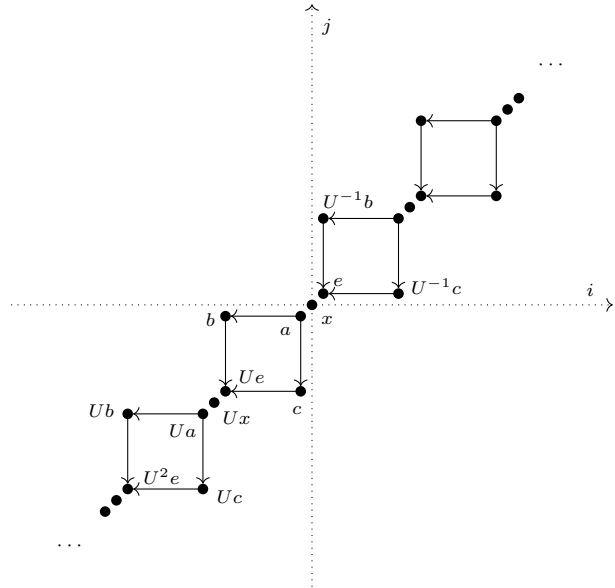

\end{example}

\subsection{The involutive concordance invariants} 
\label{sec:involutive}
We now introduce the involutive concordance invariants $\underline{V}_0$ and $\overline{V}_0$. We begin by describing an automorphism $\iota_K: CFK^{\infty}(K) \rightarrow CFK^{\infty}(K)$ for a knot $K$.

\begin{definition}
For a knot $K$ the \textit{Sarkar involution} $\sigma: CFK^{\infty}(K) \rightarrow CFK^{\infty}(K)$ is given by $$\sigma = \textrm{Id} + U^{-1}(\Phi \circ \Psi),$$ where Id is the identity map on $CFK^{\infty}(K)$, and $\Phi, \Psi: CFK^{\infty}(K) \rightarrow CFK^{\infty}(K)$ are the chain maps given as follows. Suppose $x \in CFK^{\infty}(K)$. Then, 
$$\Phi(x) = \sum_{i \mathrm{\ odd}} \partial_{i0}x$$ and
$$\Psi(x) = \sum_{j \mathrm{\ odd}} \partial_{0j}x.$$
\end{definition}

\noindent The Sarkar involution is filtered, preserves homological degree, and is an involution up to chain homotopy. We use the Sarkar map to describe the chain map $\iota_K$ on $CFK^{\infty}(K)$ as follows. For a knot $K$, the map $\iota_K: CFK^{\infty}(K) \rightarrow CFK^{\infty}(K)$ is an automorphism with the following additional properties: 
\begin{enumerate}
    \item $\iota_K$ is a skew-filtered chain map.
    \item $\iota_K$ preserves homological degree.
    \item $\iota_K^2 = \sigma$ up to chain homotopy equivalence.
\end{enumerate}

\noindent In many nice cases, although certainly not all, these properties are enough to specify $\iota_K$ up to skew-equivariant chain homotopy equivalence. Now, we review the definitions of the involutive concordance invariants. 

\begin{definition}
Let $AI_0^-$ be the mapping cone Cone$(A_0^- \xrightarrow{Q(\iota_K + \textrm{Id})} QA_0^-[-1])$. Then, the involutive concordance invariants $\underline{V}_0$ and $\overline{V}_0$ are
\[\underline{V}_0 = -\frac{1}{2}\bigg(\textrm{max }\{r: \exists x \in H_r(AI_0^-) \textrm{ s.t. } U^nx \neq 0 \textrm{ and } U^nx \notin \textrm{Im}(Q) \hspace{0.1in} \forall n\}-1\bigg),\]
\[\overline{V}_0 = -\frac{1}{2}\textrm{max }\{r: \exists x \in H_r(AI_0^-) \textrm{ s.t. } U^nx \neq 0 \hspace{0.1in} \forall n \textrm{ and } \exists m \geq 0 \textrm{ s.t. } U^mx \in \textrm{Im}(Q)\},\]
where Im$(Q)$ denotes the image of $Q$.
\end{definition}


\begin{example}
Recall that $\CFKi{(3_1)}$ has three generators which we call $a,b,$ and $c$. The boundary map $\partial$ is given by $\partial(a)=\partial(c)=0$ and $\partial(b)=Ua+c$. The map $\iota_K$ for $K = 3_1$ is a reflection over the line $i=j$. Concretely, this means that $\iota_{3_1}(a) = U^{-1}c$, $\iota_{3_1}(b) = b$, and $\iota_{3_1}(c) = Ua$. Then, the mapping cone \[AI_0^- = \textrm{Cone}(A_0^- \xrightarrow{Q(\iota_{3_1} + \textrm{Id})} QA_0^-[-1])\] looks like copies of the structure shown in Figure \ref{fig:rhtCone1}. After a change of basis, $AI_0^-$ has the simplified picture shown in Figure \ref{fig:rhtCone2}. The homology of $AI_0^-$ is $$\F_2[U] \langle [c+Qb],[Qc] \rangle.$$
\noindent $H_*(AI_0^-)$ is a module over $\F_2[U,Q]/Q^2$. So, we can form two towers from the subspaces that generate the homology as shown in Figure \ref{fig:rhtT}. We examine the homological gradings of the topmost elements of the towers. The element $c+Qb$ has homological grading -1, while $Qc$ has homological grading -2. Thus, $\underline{V}_0(3_1) = -\frac{1}{2}(-1-1) = 1$ and $\overline{V}_0(3_1) = -\frac{1}{2}(-2) = 1$.

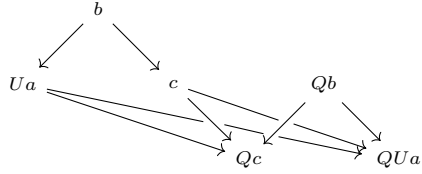
\begin{figure}
    \centering
\begin{tikzpicture}\tikzstyle{every node}=[font=\tiny]
\node(1) at (-2,0){$b$};
\node(2) at (-3,-1){$Ua$};
\node(3) at (-1, -1){$c$};

\path[->](1)edge(2);
\path[->](1)edge(3);

\node(4) at (1,-1){$Qb$};
\node(5) at (2,-2){$QUa$};
\node(6) at (0,-2){$Qc$};

\path[->](4)edge(5);
\path[->](4)edge(0.2,-1.8);

\path[->](3)edge(-0.23,-1.75);

\path[->](2)edge(6);

\draw (3) -- (0.4,-1.48);
\draw[->](0.6,-1.53) -- (5);

\draw (2) -- (-0.6,-1.5);
\draw (-0.32,-1.54) -- (0.2,-1.66);
\draw[->](0.4,-1.7) -- (1.51,-1.93);
\end{tikzpicture}
    \caption{$AI_0^-$ for the right-handed trefoil knot consists of all positive $U$-powers of the structure shown.}
    \label{fig:rhtCone1}
\end{figure}

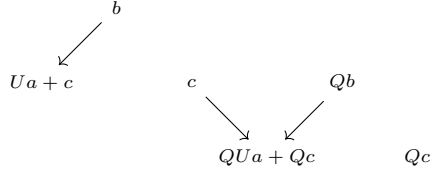
\begin{figure}
    \centering
\begin{tikzpicture}\tikzstyle{every node}=[font=\tiny]
\node(1) at (-2,0){$b$};
\node(2) at (-3,-1){$Ua+c$};
\node(3) at (-1,-1){$c$};
\path[->](1)edge(2);
\node(4) at (1,-1){$Qb$};
\node(5) at (0,-2){$QUa+Qc$};
\node(6) at (2,-2){$Qc$};
\path[->](4)edge(5);
\path[->](3)edge(5);
\end{tikzpicture}
    \caption{$AI_0^-$ for the right-handed trefoil knot shown after applying a change of basis to the structure in Figure \ref{fig:rhtCone1}.}
    \label{fig:rhtCone2}
\end{figure}

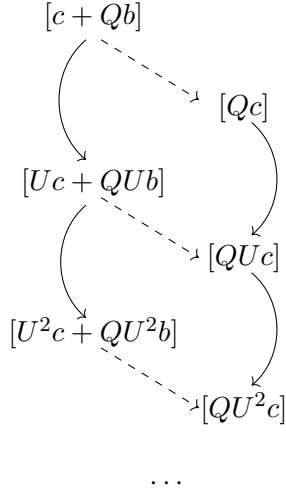
\begin{figure}
    \centering
    \begin{tikzpicture}
    \node(1)at(0,0.2){$[c+Qb]$};
    \path[->][bend right = 50](-.1,-.1)edge(-.1,-1.7);
    \path[->][dashed](.1,-.1)edge(1.4,-0.9);
    \node(2)at(0,-2){$[Uc+QUb]$};
    \path[->][bend right = 50](-.1,-2.3)edge(-.1,-3.8);
    \path[->][dashed](.1,-2.2)edge(1.4,-3);
    \node(3)at(0,-4){$[U^2c+QU^2b]$};
    \path[->][dashed](.1,-4.2)edge(1.4,-5);
    \node(5)at(2,-1){$[Qc]$};
    \path[->][bend left = 50](2.1,-1.2)edge(2.1,-2.7);
    \node(6)at(2,-3){$[QUc]$};
    \path[->][bend left = 50](2.1,-3.2)edge(2.1,-4.7);
    \node(7)at(2,-5){$[QU^2c]$};
    \node(8)at(1,-6){$\bf \cdots$};
    \end{tikzpicture}
    \caption{$H_*(AI_0^-)$ for the right-handed trefoil knot can be described by two linked towers. The tower on the left contains no elements in the image of $Q$, while the tower on the right contains only elements that are in the image of $Q$. Each tower is organized by increasing powers of $U$. Curved lines denote multiplication by $U$, while dashed lines denote application of $Q$.}
    \label{fig:rhtT}
\end{figure}
\end{example}


\begin{example}
The involutive concordance invariants for the left-handed trefoil knot can be calculated in a similar way. Recall that $\CFKi{(\overline{3_1})}$ has three generators which we call $a,b,$ and $c$. The boundary map $\partial$ is given by $\partial(a)= b, \partial(b)=0$, and $\partial(c)=Ub$. As with the right handed trefoil, $\iota(\overline{3}_1)$ is a reflection over the line $i=j$. The mapping cone \[AI_0^- = \textrm{Cone}(A_0^- \xrightarrow{Q(\iota_{\overline{3}_1} + \textrm{Id})} QA_0^-[-1])\] has the form shown in Figure \ref{fig:lhtCone1}. The result of changing the basis is shown in Figure \ref{fig:lhtCone2}. The homology of $AI_0^-$ is $$\F_2[U] \langle [Ua+c],[QUa+Qc] \rangle \cup [b] \cup [Qb].$$
\noindent We can form two towers from the subspaces that generate the homology as shown in Figure \ref{fig:lhtT}. We examine the homological gradings of the topmost elements of the towers. The element $Ua+c$ has homological grading 1, while $b$ has homological grading 2. Thus, $\underline{V}_0(\overline{3_1}) = -\frac{1}{2}(1-1) = 0$ and $\overline{V}_0(\overline{3_1}) = -\frac{1}{2}(2) = -1$.

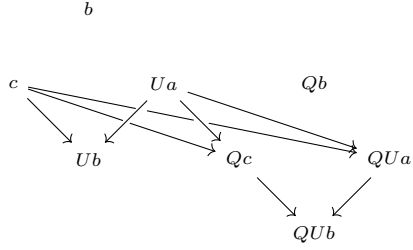
\begin{figure}
    \centering
    \begin{tikzpicture}\tikzstyle{every node}=[font=\tiny]
    \node(0) at (-1, 0){$b$};
    \node(1) at (-1,-2){$Ub$};
    \node(2) at (-2,-1){$c$};
    \node(3) at (-0, -1){$Ua$};
    
    \path[->](2)edge(1);
    \path[->](3)edge(1);
    
    \node(7) at (2, -1){$Qb$};
    \node(4) at (2,-3){$QUb$};
    \node(5) at (1,-2){$Qc$};
    \node(6) at (3, -2){$QUa$};
    
    \path[->](5)edge(4);
    \path[->](6)edge(4);
    
    \path[->](3)edge(6);
    \path[->](3)edge(5);
    
    
    \draw (2) -- (-0.54,-1.485);
    \draw[->](-0.45,-1.52) -- (5);
    
    \draw (2) -- (-0.37,-1.325);
    \draw(-0.22,-1.355) -- (0.4,-1.48);
    \draw[->](0.6,-1.52) -- (6);
    \end{tikzpicture}
    \caption{The mapping cone for the left handed trefoil.}
    \label{fig:lhtCone1}
\end{figure}

\begin{figure}
    \centering
\begin{tikzpicture}\tikzstyle{every node}=[font=\tiny]
    \node(0) at (-1, 1){$b$};
    \node(1) at (-1,-1){$Ub$};
    \node(2) at (-2,0){$Ua+c$};
    \node(3) at (0, 0){$Ua$};
    \path[->](3)edge(1);
    \node(7) at (3, 0){$Qb$};
    \node(4) at (3,-2){$QUb$};
    \node(5) at (4,-1){$Qc$};
    \node(6) at (2, -1){$QUa+Qc$};
    \path[->](5)edge(4);
    \path[->](3)edge(6);
    \end{tikzpicture}
    \caption{The simplified mapping cone for the left handed trefoil.}
    \label{fig:lhtCone2}
\end{figure}
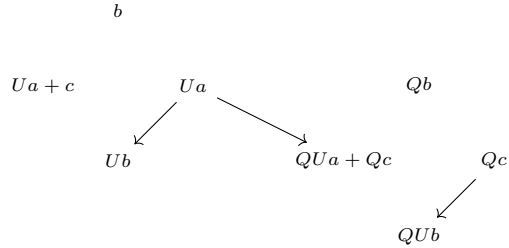
 
 \begin{figure}
     \centering
    \begin{tikzpicture}
    \node(1)at(0,0.2){$[Ua+c]$};
    \path[->][bend right = 50](-.1,-.1)edge(-.1,-1.7);
    \path[->][dashed](.1,-.1)edge(1.4,-0.9);
    \node(2)at(0,-2){$[U^2a+Uc]$};
    \path[->][bend right = 50](-.1,-2.3)edge(-.1,-3.8);
    \path[->][dashed](.1,-2.2)edge(1.4,-3);
    \node(3)at(0,-4){$[U^3a+U^2c]$};
    \path[->][dashed](.1,-4.2)edge(1.4,-5);
    \node(4)at(2.5,1){$[b]$};
    \path[->][bend left = 50](2.6,0.7)edge(2.6,-0.7);
    \node(5)at(2.5,-1){$[QUa+Qc]$};
    \path[->][bend left = 50](2.6,-1.2)edge(2.6,-2.7);
    \node(6)at(2.7,-3){$[QU^2a+QUc]$};
    \path[->][bend left = 50](2.6,-3.2)edge(2.6,-4.7);
    \node(7)at(2.8,-5){$[QU^3a+QU^2c]$};
    \node(9)at(4.5,0){$[Qb]$};
    \node(8)at(1,-6){$\bf \cdots$};
    \end{tikzpicture}
     \caption{$H_*(AI_0^-)$ for the left handed trefoil knot can be described by two linked towers and the stand-alone subspace $[Qb]$. Note that $QUb$ is in the image of the differential.}
     \label{fig:lhtT}
 \end{figure}
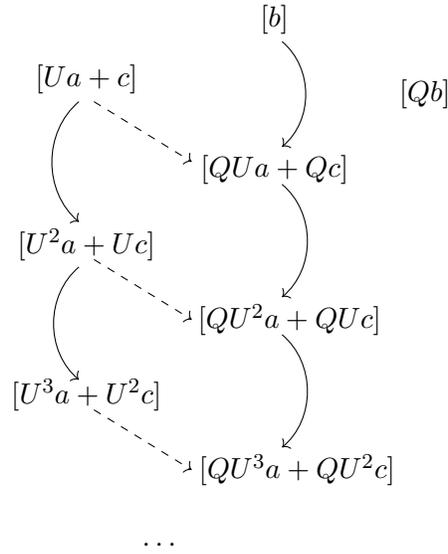
\end{example}


\begin{example}
We also find the involutive concordance invariants for the figure-eight knot from Example \ref{ex:figeight}. The map $\iota_K$ for $K = 4_1$ is given by $\iota_{4_1}(a) = a+x$, $\iota_{4_1}(b) = c, \iota_{4_1}(c) = b$, $\iota_{4_1}(e) = e$, and $\iota_{4_1}(x) = e+x$. Then, the mapping cone \[AI_0^- = \textrm{Cone}(A_0^- \xrightarrow{Q(\iota_{4_1} + \textrm{Id})} QA_0^-[-1])\] is represented by the structure shown in Figure \ref{fig:figeightCone1}. After a change of basis, $AI_0^-$ has the simplified picture shown in Figure \ref{fig:figeightCone2}. The homology of $AI_0^-$ is $$\F_2[U] \langle [Ux+Qc],[Qx] \rangle \cup [e].$$
\noindent Thus, we can form two towers from the subspaces that generate the homology as shown in Figure \ref{fig:figeightT}. We examine the homological gradings of the topmost elements of the towers. The element $Ux + Qc$ has homological grading -1, while $Qx$ has homological grading 0. Thus, $\underline{V}_0(4_1) = -\frac{1}{2}(-1-1) = 1$ and $\overline{V}_0(3_1) = -\frac{1}{2}(0) = 0$.

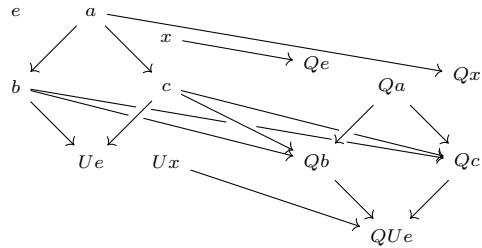
\begin{figure}
    \centering
    \begin{tikzpicture}\tikzstyle{every node}=[font=\tiny]
    \node(1) at (-2,0){$a$};
    \node(-1) at (-3,0){$e$};
    \node(0) at (-1,-0.35){$x$};
    \node(3) at (-3,-1){$b$};
    \node(4) at (-1,-1){$c$};
    \node(2) at (-2,-2){$Ue$};
    \node(11) at (-1,-2){$Ux$};
    
    \path[->](1)edge(4);
    \path[->](1)edge(3);
    
    \path[->](3)edge(2);
    \path[->](4)edge(2);
    \node(7) at (2,-1){$Qa$};
    \node(5) at (1,-0.7){$Qe$};
    \node(6) at (3,-0.85){$Qx$};
    \node(9) at (1,-2){$Qb$};
    \node(10) at (3,-2){$Qc$};
    \node(8) at (2,-3){$QUe$};
    
    \path[->](7)edge(9);
    \path[->](7)edge(10);
    
    \path[->](9)edge(8);
    \path[->](10)edge(8);
    
    \path[->](11)edge(8);
    
    \path[->](0)edge(5);
    \path[->](1)edge(6);
    
    \draw (4) -- (1.3,-1.58);
    \draw[->] (1.5,-1.63) -- (10);
    
    \path[->](4)edge(9);
    
    \draw (3) -- (-1.5,-1.38);
    \draw[->] (-1.3,-1.42) -- (9);
    
    \draw (3) -- (-1.4,-1.27);
    \draw (-1.2,-1.3) -- (-0.2,-1.47);
    \draw (0.2,-1.53) -- (1.2,-1.7);
    \draw[->] (1.4,-1.735) -- (10);
    
    \end{tikzpicture}
    \caption{$AI_0^-$ for the figure eight knot consists of the structure shown, and all positive $U$-translates of the generators $a, b, c, x, Ux, Ue, Qa, Qb, Qc, Qx,$ and $QUe$.}
    \label{fig:figeightCone1}
\end{figure}

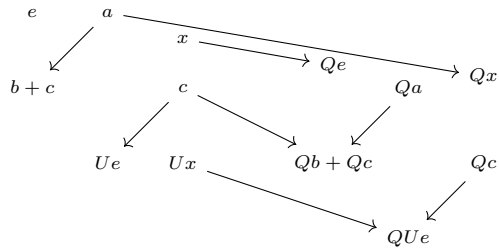
\begin{figure}
    \centering
    \begin{tikzpicture}\tikzstyle{every node}=[font=\tiny]
    \node(-1) at (-3,0){$e$};
    \node(0) at (-1,-0.35){$x$};
    \node(1) at (-2,0){$a$};
    \node(2) at (-2,-2){$Ue$};
    \node(3) at (-3,-1){$b+c$};
    \node(4) at (-1,-1){$c$};
    \node(11) at (-1,-2){$Ux$};
    \path[->](1)edge(3);
    \path[->](4)edge(2);
    \node(5) at (1,-0.7){$Qe$};
    \node(6) at (3,-0.85){$Qx$};
    \node(7) at (2,-1){$Qa$};
    \node(8) at (2,-3){$QUe$};
    \node(9) at (1,-2){$Qb+Qc$};
    \node(10) at (3,-2){$Qc$};
    \path[->](7)edge(9);
    \path[->](10)edge(8);
    \path[->](0)edge(5);
    \path[->](1)edge(6);
    \path[->](4)edge(9);
    \path[->](11)edge(8);
    \end{tikzpicture}
    \caption{$AI_0^-$ for the figure eight knot shown after applying a change of basis to the structure in Figure \ref{fig:figeightCone1}.}
    \label{fig:figeightCone2}
\end{figure}

\begin{figure}
    \centering
    \begin{tikzpicture}
    \node(0)at(4,2){$[e]$};
    \node(1)at(0,0.2){$[Ux + Qc]$};
    \path[->][bend right = 50](-.1,-.1)edge(-.1,-1.7);
    \path[->][dashed](.1,-.1)edge(1.4,-0.9);
    \node(2)at(0,-2){$[U^2x + QUc]$};
    \path[->][bend right = 50](-.1,-2.3)edge(-.1,-3.8);
    \path[->][dashed](.1,-2.2)edge(1.4,-3);
    \node(3)at(0,-4){$[U^3x + QU^2c]$};
    \path[->][dashed](.1,-4.2)edge(1.4,-5);
    \node(4)at(2,1){$[Qx]$};
    \path[->][bend left = 50](2.1,0.7)edge(2.1,-0.7);
    \node(5)at(2,-1){$[QUx]$};
    \path[->][bend left = 50](2.1,-1.2)edge(2.1,-2.7);
    \node(6)at(2,-3){$[QU^2x]$};
    \path[->][bend left = 50](2.1,-3.2)edge(2.1,-4.7);
    \node(7)at(2,-5){$[QU^3x]$};
    \node(8)at(1,-6){$\bf \cdots$};
    \end{tikzpicture}
    \caption{$H_*(AI_0^-)$ for the figure eight knot can be described by two linked towers and the stand-alone subspace $[e]$.}
    \label{fig:figeightT}
\end{figure}
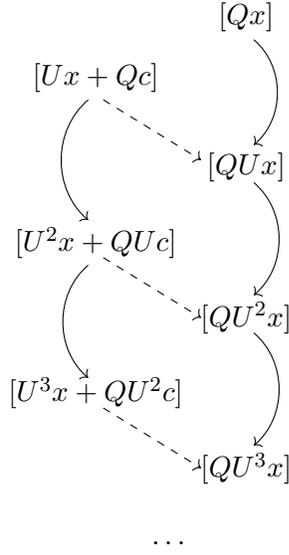

\end{example}


\section{Heegaard Diagrams for $(1,1)$-Knots}\label{sec:1,1}

\begin{figure}
    \centering
    \includegraphics[width=0.4\textwidth]{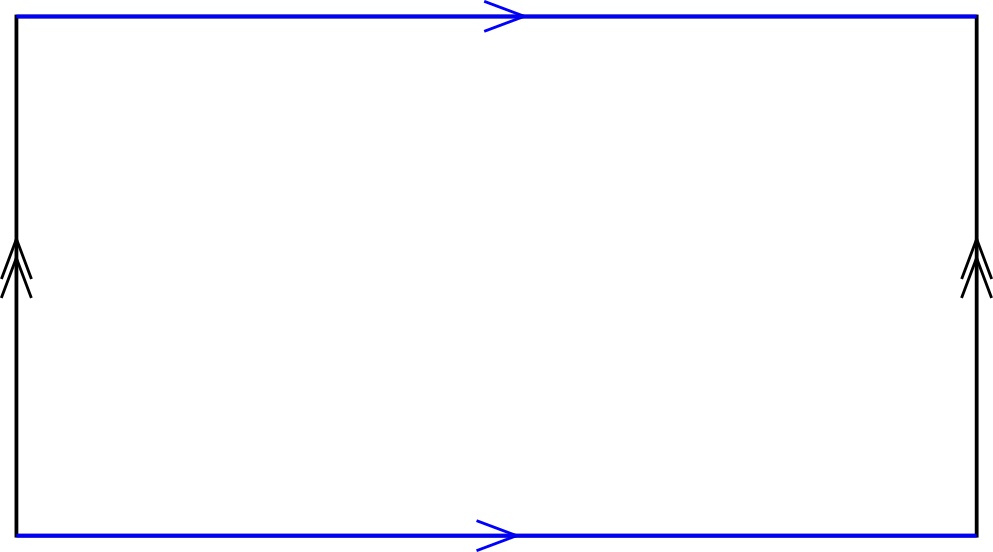}
    \caption{The torus can be represented as a rectangle with opposite sides identified.}
    \label{fig:rectangle}
\end{figure}

We can compute $\CFKi{(K)}$ for a (1,1)-knot $K$ using information from a representation of $K$ on the torus, called a Heegaard diagram. Heegaard diagrams for $(1,1)$-knots are uniquely determined by a 4-tuple of integers, described by \cite{Racz} as follows:

\begin{definition}
For a (1,1)-knot, there exists a (nonunique) parameterization $(k,r,c,s)$ that describes how the knot lies on the torus, where:
\begin{enumerate}
    \item There are closed curves $\alpha$ and $\beta$ such that there are $2k+1$ intersections between $\alpha$ and $\beta$, labeled $x_0,x_1,...,x_k,x_{-k},...,x_{-1},x_{0}$. Visualizing the torus as a rectangle with identified sides as in Figure \ref{fig:rectangle}, $\beta$ consists of the top, or equivalently the bottom, of the rectangle, and $\alpha$ is the union of the following arcs. In the following, the words `left' and `right' are with respect to the orientation of the curve.
    \item There are $r$ loops connecting $x_{c-i}$ to $x_{c+i}$ on the left side of $\beta$, for $k-r<i \leq k$. The centermost loop contains the basepoint $w$.
    \item There are $r$ loops connecting $x_{-(c-i)}$ to $x_{-(c+i)}$ on the right side of $\beta$, for $k-r<i \leq k$. The centermost loop contains the basepoint $z$.
    \item There are $|s|$ bridges connecting the left side of $\beta$ at $x_i$ and the right side of $\beta$ at $x_j$, where:
    \begin{itemize}
        \item $c-k+r \leq i < c-k+r+|s|$,
        \item $-(c-k+r+|s|)<j \leq -(c-k+r)$, and
        \item $i-j=2(c-k+r)+|s|-1$.
    \end{itemize}
    \item There are $t = 2(k-r) + 1 - |s|$ bridges connecting the left side of $\beta$ at $x_i$ and the right side of $\beta$ at $x_j$ where:
    \begin{itemize}
        \item $c+k-r-t<i \leq c+k-r$,
        \item $-(c+k-r) \leq j \leq -(c+k-r-t)$, and
        \item $i-j=2(c+k-r)-t+1.$
    \end{itemize}
\end{enumerate}
\end{definition}

\begin{figure}
    \centering
    \includegraphics[width=0.6\textwidth]{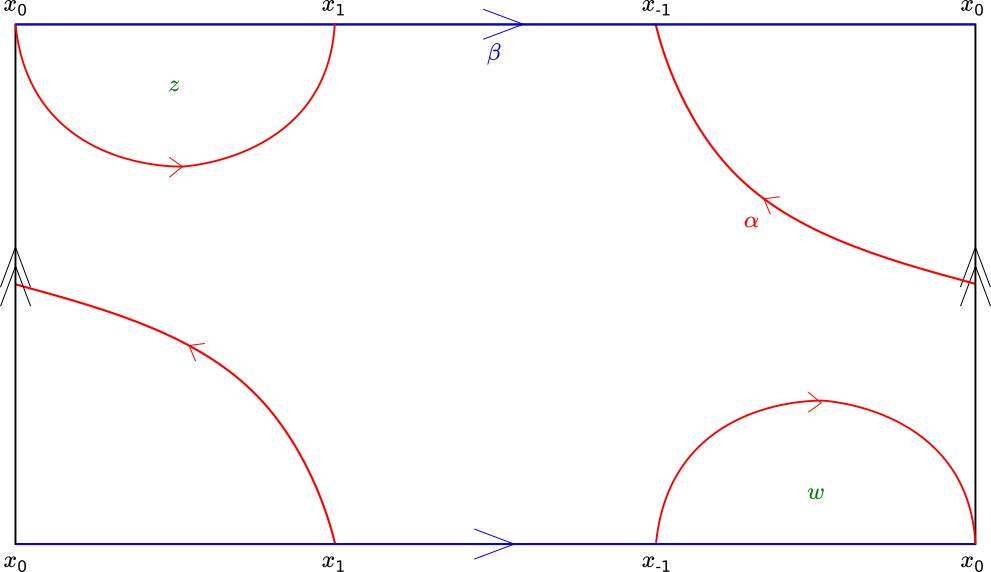}
    \caption{The Heegaard diagram for the right-handed trefoil knot.}
    \label{fig:heegaard_trefoil}
\end{figure}

\noindent There is also a \textit{Rasmussen parameterization} for $(1,1)$-knots which is given by a different 4-tuple \cite{rasmussen}. A Heegaard diagram for the right-handed trefoil knot $3_1$ is shown in Figure \ref{fig:heegaard_trefoil}.

\begin{figure}
    \centering
    \includegraphics[width=0.6\textwidth]{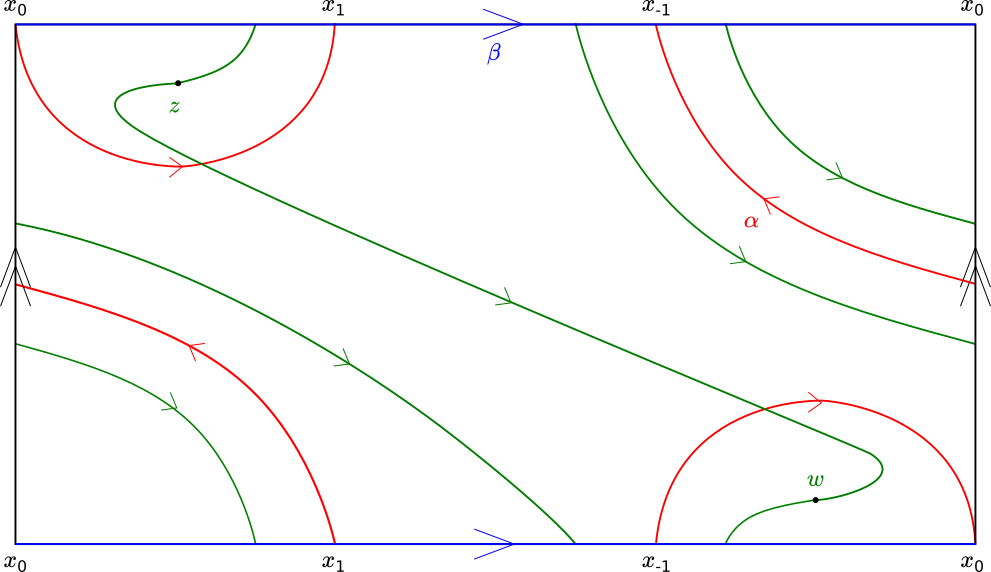}
    \caption{The Heegaard diagram for the right-handed trefoil knot, with the knot shown in green.}
    \label{fig:heegaard_with_trefoil}
\end{figure}

A Heegaard diagram for a (1,1)-knot $K$ determines how $K$ lies on a torus in the following way. The understrand of $K$ is drawn from $w$ to $z$ without intersecting $\alpha$. Then, the overstrand of $K$ is drawn from $z$ to $w$ without intersecting $\beta$. Figure \ref{fig:heegaard_with_trefoil} depicts the right-handed trefoil knot on the torus.

\begin{figure}
    \centering
    \includegraphics[width=.15\linewidth]{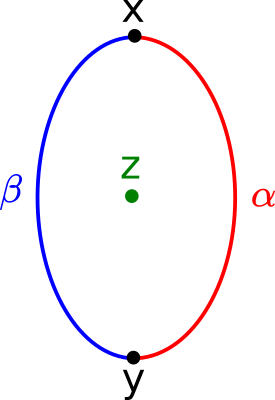} \vspace{1in}
    \caption{A bigon of this form corresponds to an appearance of the element $x$ in the boundary of $y$. Furthermore, since the bigon contains a single copy of $z$, $A(x) - A(y) = 1$ and $M(x) - M(y) = 1$.}
    \label{fig:bigon}
\end{figure}

We may extract the generators for $\CFKi{(K)}$ from the Heegaard diagram, and identify bigons in the diagram, which determine the boundary map. A \textit{bigon} in the $(1,1)$-diagram for a knot $K$ is a disk on the torus $\Sigma$ whose boundary consists of one segment from each of the curves $\alpha$ and $\beta$ such that $\alpha$ and $\beta$ intersect exactly twice. The disk must be convex at the intersection points, that is, it must occupy one of the four regions of $\Sigma - \alpha - \beta$ which meet at that corner. Examples of bigons are shown in Figures \ref{fig:bigon}, \ref{fig:Bigon}, and \ref{fig:Bigon2}. Bigons determine the boundary map on $\CFKi{(K)}$ in the following way. Given a bigon with intersection points $x$ and $y$ between the curves $\alpha$ and $\beta$, we orient the bigon so that the part of the curve $\beta$ on the boundary is on the left. Assume that with this orientation, $x$ is on the top and $y$ is on the bottom. Suppose the bigon contains $m$ copies of the point $z$ and $n$ copies of the point $w$. Then, the bigon corresponds to an appearance of the element $U^nx$ in the boundary of $y$. Moreover, the difference in Alexander gradings is $A(x) - A(y) = m - n$. The difference in homological gradings is $M(x) - M(y) = 1 - 2n$. By taking the sum over all bigons in a Heegaard diagram for a knot $K$, the differential $\partial$ on $\CFKi{(K)}$ is determined. An example of this process is carried out in Section \ref{sub:10_161cfk}.

\section{Example Computations of the Involutive Concordance Invariants}
\label{sec:examples}
We computed the involutive concordance invariants for each of the $(1,1)$-knots included in Table \ref{tab:invol} in Section \ref{sec:results}. We chose to consider these 10- and 11-crossing (1,1)-knots based on the following classification.

The 10-crossing (1,1)-knots were classified by Morimoto, Sakuma, and Yokota in \cite{MSYtunnel}, and the resulting list appears in \cite[Table 1]{GMM}. The knots $10n_{125}$, $10n_{126}$, $10n_{127}$, $10n_{129}$, $10n_{130}$, $10n_{131}$, $10n_{133}$, $10n_{134}$, $10n_{135}$, $10n_{137}$, and $10n_{138}$ are \emph{thin} in the terminology of Heegaard Floer homology. The knot $10n_{124}$ is an \emph{$L$-space knot} in the terminology of Heegaard Floer homology. The involutive concordance invariants for thin knots and $L$-space knots have already been determined in \cite[Section 8]{HM:involutive}. The remaining knots, $10n_{128}$, $10n_{132}$, $10n_{136}$, $10n_{139}$, $10n_{145}$, and $10n_{161}$, are of interest.

 We grouped 11-crossing knots based on whether or not they are \emph{Montesinos}. Montesinos knots are knots composed of rational tangles. Racz classified which 11-crossing non-Montesinos knots are (1,1) in \cite[Section 3]{Racz}. These knots are $11n96$, $11n111$, and $11n135$. These knots are neither thin nor $L$-space. Klimenko and Sakuma classified which 11-crossing Montesinos knots are (1,1) by showing that they are exactly the 11-crossing Montesinos knots with tunnel number one \cite[Corollary C]{KSmontesinos}. Castellano-Mac\'{i}as and Owad list all 11-crossing knots with tunnel number one in \cite[Appendix A]{tunnelNumber}. Only the non-alternating knots are interesting, since alternating knots are thin. The thin, non-alternating knots on this list are $11n_{1}$, $11n_{2}$, $11n_{3}$, $11n_{13}$, $11n_{14}$, $11n_{15}$, $11n_{16}$, $11n_{17}$, $11n_{18}$, $11n_{28}$, $11n_{29}$, $11n_{30}$, $11n_{51}$, $11n_{52}$, $11n_{53}$, $11n_{54}$, $11n_{55}$, $11n_{56}$, $11n_{58}$, $11n_{59}$, $11n_{60}$, $11n_{62}$, $11n_{63}$, and $11n_{64}$. There are no $L$-space knots on the list. The knots $11n_{143}$ and $11n_{145}$ are non-Montesinos. So, the remaining knots of interest are $11n_{12}$, $11n_{19}$, $11n_{20}$, $11n_{38}$, $11n_{57}$, $11n_{61}$, $11n_{70}$, $11n_{79}$, $11n_{102}$, and $11n_{104}$. Thus, we have the final list of 10- and 11-crossing (1,1)-knots for which we calculated the involutive concordance invariants.

The first part of calculating the involutive concordance invariants for the selected (1,1)-knots involved finding a model for $\CFKi{(K)}$ for each knot $K$. We used one of two strategies to complete this step. For the knots $10_{128}, 10_{132}, 10_{136}, 10_{139}, 10_{145}, 11n_{12}, 11n_{19}, 11n_{57}, 11n_{70}$, and $11n_{79}$ we used the Heegaard Floer knot homology of each knot given in Section 3 of \cite{hfkhat10cross} to find $\CFKi{(K)}$; in these cases, this was enough information to specify the full complex up to chain homotopy. For the remaining knots, we used the information in the table in Section 3.7.3 of \cite{Racz} to draw the $(1,1)$ Heegaard diagram for each knot $K$, and use this diagram to compute $\CFKi{(K)}.$

\subsection{Finding the involutive concordance invariants for $11n_{57}$.}\label{sub:11n57ex} As an example of the first strategy, we consider the computation of $\underline{V}_0(11n_{57})$ and $\overline{V}_0(11n_{57})$. The first step is to find $\CFKi{(11n_{57})}$ using the Heegaard Floer knot homology of the knot.

\subsubsection{Finding $\CFKi{(11n_{57})}$ using Heegaard Floer knot homology} 
\label{sub:fromHFKhat}

\begin{figure}
    \centering
    \begin{tikzpicture}\tikzstyle{every node}=[font=\tiny] 
    \path[->][dotted](0,-5)edge(0,5.3);
    \path[->][dotted](-4,0)edge(4,0);
    \node() at (.2,5){$j$};
    \node() at (3.7,.2){$i$};
    
	\fill(0,4)circle [radius=2pt];
    \node(8) at (0.15,4.2){$a$};     
    
    \fill(-0.2,3)circle[radius=2pt];
    \node(7) at (-0.2,3.2){$b$};
    \fill(0,3)circle [radius=2pt];
    \node(6) at (0.1,3.2){$c$};
    \fill(0.2,3)circle [radius=2pt];
    \node(5) at (0.4,3.2){$d$};    
    
	\fill(-0.1,2)circle[radius=2pt];
    \node(4) at (-0.2,2.2){$e$};
    \fill(0.1,2)circle [radius=2pt];
    \node(3) at (0.2,2.2){$f$};    
    
	\fill(0,1)circle [radius=2pt];
    \node(2) at (0.1,1.2){$g$};    
    
    \fill(-0.2,0)circle[radius=2pt];
    \node(-1) at (-0.2,0.2){$h$};
    \fill(0,0)circle [radius=2pt];
    \node(0) at (0.1,0.2){$i$};
    \fill(0.2,0)circle [radius=2pt];
    \node(1) at (0.4,0.2){$j$};
    
    \fill(0,-1)circle [radius=2pt];
    \node(-2) at (0.1,-1.2){$k$}; 
    
    \fill(-0.1,-2)circle[radius=2pt];
    \node(-4) at (-0.2,-2.2){$l$};
    \fill(0.1,-2)circle [radius=2pt];
    \node(-3) at (0.2,-2.2){$m$};
    
    \fill(-0.2,-3)circle[radius=2pt];
    \node(-7) at (-0.2,-3.2){$n$};
    \fill(0,-3)circle [radius=2pt];
    \node(-6) at (0.1,-3.2){$o$};
    \fill(0.2,-3)circle [radius=2pt];
    \node(-5) at (0.4,-3.2){$p$}; 
    
    \fill(0,-4)circle [radius=2pt];
    \node(-8) at (0.15,-4.2){$q$}; 
    \end{tikzpicture}
    \caption{Arrangement of generators of $\CFKi{(11n_{57})}$. The homological grading of $a$ is 1; $b$, $c$, and $d$ are in homological grading 0; $e,g,$ and $g$ are in homological grading $-1$; $h,i$, and $j$ are in homological grading $-2$; $k$ is in homological grading $-3$; $l$ and $m$ are in homological grading $-5$; $n,o,$ and $p$ are in homological grading $-6$; $q$ has homological grading $-7$.}
    \label{fig:cfk57_1}
\end{figure}
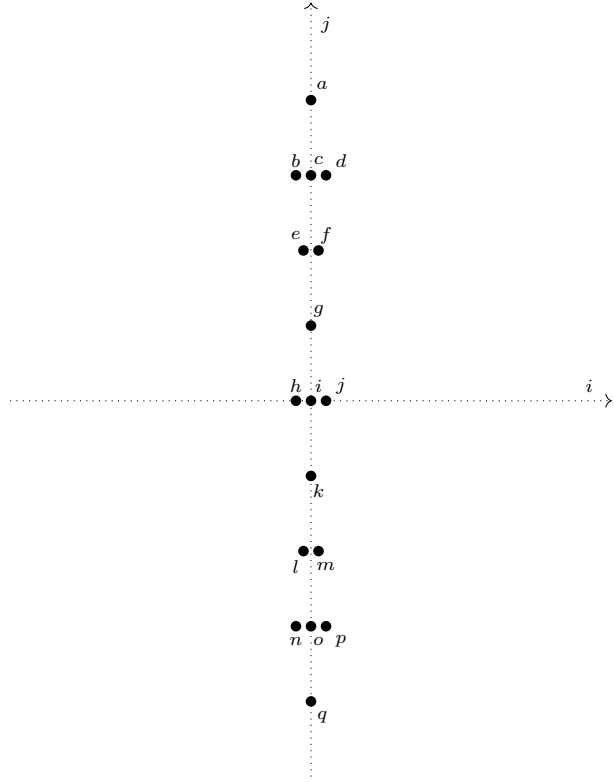

\begin{figure}
    \centering
    \begin{tikzpicture}\tikzstyle{every node}=[font=\tiny] 
    \path[->][dotted](0,-4)edge(0,4.3);
    \path[->][dotted](-6,0)edge(6,0);
    \node() at (.2,4){$j$};
    \node() at (5.7,.2){$i$};
    
	\fill(-1,3)circle [radius=2pt];
    \node at (-1.15,3.2){$Ua$};
    \path[->](-1,3)edge(-1,2.2);  
    
    \fill(-1,2.1)circle[radius=2pt];
    \node at (-1.35,2.15){$Ub$};
    \fill(-1,1.9)circle [radius=2pt];
    \node at (-1.35,1.85){$Uc$};
    \fill(0,3)circle [radius=2pt];
    \node at (0.2,3.2){$d$};
    \path[->](0,3)edge(-0.9,3);
    \path[->](0,3)edge(0,2.2);
    
	\fill(-0,2.1)circle[radius=2pt];
    \node at (0.2,2.1){$f$};
    \path[->](0,2.1)edge(-0.9,2.1);
    \fill(0,1.9)circle [radius=2pt];
    \node at (0.2,1.8){$e$};
    \path[->](0,1.9)edge(-0.9,1.9);
    \path[->](0,1.9)edge(0,0.1); 
    
	\fill(0.2,1.2)circle [radius=2pt];
    \node at (0.3,1.4){$g$};
    \path[->](0.2,1.2)edge(0.2,0.3);  
    
    \fill(0,0)circle[radius=2pt];
    \node at (-0.2,0.2){$h$};
    \fill(0.2,0.2)circle [radius=2pt];
    \node at (0.4,0.4){$i$};
    \fill(1.2,1.2)circle [radius=2pt];
    \node at (1.65,1.4){$U^{-1}j$};
    \path[->](1.2,1.2)edge(1.2,0.3);
    \path[->](1.2,1.2)edge(0.3,1.2);
    
    \fill(1.2,0.2)circle [radius=2pt];
    \node at (1.65,0.5){$U^{-1}k$};
    \path[->](1.2,0.2)edge(0.3,0.2);
    
    \fill(2.4,0)circle[radius=2pt];
    \node at (2.8,0.2){$U^{-2}l$};
    \path[->](2.4,0)edge(2.4,-0.9);
    \fill(2.2,0)circle [radius=2pt];
    \node at (1.7,-0.2){$U^{-2}m$};
    \path[->](2.2,0)edge(0.1,0);
    \path[->](2.2,0)edge(2.2,-0.9);
    
    \fill(2.4,-1)circle[radius=2pt];
    \node at (2.7,-1.2){$U^{-2}n$};
    \fill(2.2,-1)circle [radius=2pt];
    \node at (1.7,-1.1){$U^{-2}o$};
    \fill(3.4,0)circle [radius=2pt];
    \node at (3.8,0.2){$U^{-3}p$}; 
    \path[->](3.4,0)edge(2.5,0);
    \path[->](3.4,0)edge(3.4,-0.9);
    \fill(3.4,-1)circle [radius=2pt];
    \node at (3.8,-1.2){$U^{-3}q$};
    \path[->](3.4,-1)edge(2.5,-1);
    \end{tikzpicture}
    \caption{$\CFKi{(11n_{57})}$ decomposes as an $\F_2[U,U^{-1}]$ complex into the direct sum of a staircase and three boxes.}
    \label{fig:cfk57_final}
\end{figure}
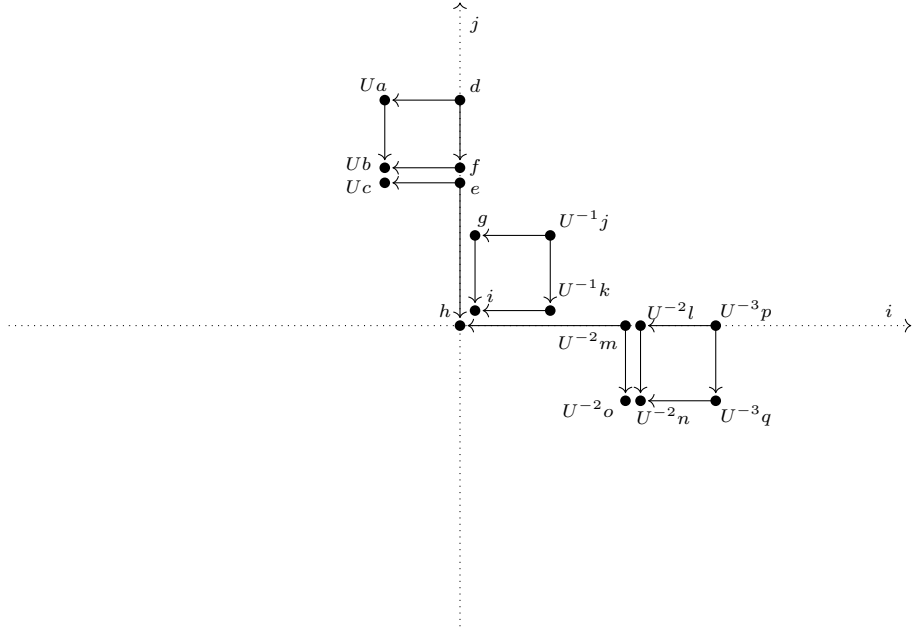

We detail how to find $\CFKi{(11n_{57})}$, which is shown in Figure \ref{fig:cfk57_final}. The Poincar\'{e} polynomial of the knot Floer homology of $11n_{57}$ is given by 
\[\HFKhat{(S^3,11n_{57})} = q^{-7}t^{-4} + 3q^{-6}t^{-3} + 2q^{-5}t^{-2} + q^{-3}t^{-1} + 3q^{-2} + q^{-1}t + 2q^{-1}t^2 + 3t^3 + qt^4, \]

\noindent where an entry $q^mt^n$ in the sum denotes a one-dimensional summand in the homology in homological grading $m$ and Alexander grading $n$. So, the generators of $\CFKi{(11n_{57})}$ are arranged as in Figure \ref{fig:cfk57_1}.

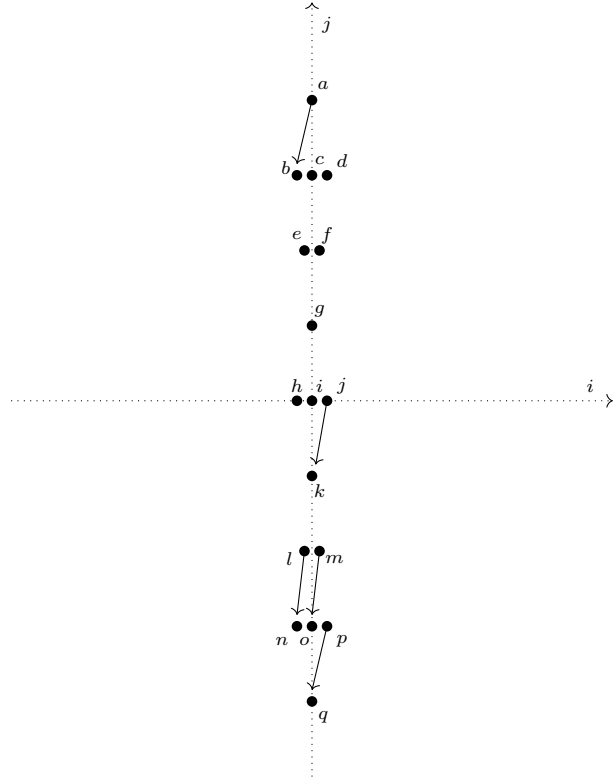
\begin{figure}
    \centering
    \begin{tikzpicture}\tikzstyle{every node}=[font=\tiny] 
    \path[->][dotted](0,-5)edge(0,5.3);
    \path[->][dotted](-4,0)edge(4,0);
    \node() at (.2,5){$j$};
    \node() at (3.7,.2){$i$};
    
	\fill(0,4)circle [radius=2pt];
    \node(8) at (0.15,4.2){$a$};
    \path[->](0,4)edge(-0.2,3.15);  
    
    \fill(-0.2,3)circle[radius=2pt];
    \node(7) at (-0.35,3.1){$b$};
    \fill(0,3)circle [radius=2pt];
    \node(6) at (0.1,3.2){$c$};
    \fill(0.2,3)circle [radius=2pt];
    \node(5) at (0.4,3.2){$d$};    
    
	\fill(-0.1,2)circle[radius=2pt];
    \node(4) at (-0.2,2.2){$e$};
    \fill(0.1,2)circle [radius=2pt];
    \node(3) at (0.2,2.2){$f$};    
    
	\fill(0,1)circle [radius=2pt];
    \node(2) at (0.1,1.2){$g$};    
    
    \fill(-0.2,0)circle[radius=2pt];
    \node(-1) at (-0.2,0.2){$h$};
    \fill(0,0)circle [radius=2pt];
    \node(0) at (0.1,0.2){$i$};
    \fill(0.2,0)circle [radius=2pt];
    \node(1) at (0.4,0.2){$j$};
    \path[->](0.2,0)edge(0.05,-0.85);
    
    \fill(0,-1)circle [radius=2pt];
    \node(-2) at (0.1,-1.2){$k$}; 
    
    \fill(-0.1,-2)circle[radius=2pt];
    \node(-4) at (-0.3,-2.1){$l$};
    \path[->](-0.1,-2)edge(-0.2,-2.85);
    \fill(0.1,-2)circle [radius=2pt];
    \node(-3) at (0.3,-2.1){$m$};
    \path[->](0.1,-2)edge(0,-2.85);
    
    \fill(-0.2,-3)circle[radius=2pt];
    \node(-7) at (-0.4,-3.2){$n$};
    \fill(0,-3)circle [radius=2pt];
    \node(-6) at (-0.1,-3.2){$o$};
    \fill(0.2,-3)circle [radius=2pt];
    \node(-5) at (0.4,-3.2){$p$};
    \path[->](0.2,-3)edge(0,-3.85);
    
    \fill(0,-4)circle [radius=2pt];
    \node(-8) at (0.15,-4.2){$q$}; 
    \end{tikzpicture}
    \caption{An incomplete picture for $\CFKi{(11n_{57})}$ with part of $\partial_{\textrm{vert}}$ denoted by the arrows. Here $\partial(c) = 0$ and the vertical differentials of the remaining unpaired elements are not yet determined.}
    \label{fig:cfk57_2}
\end{figure}

We next inspect the vertical differential. By Lemma 2.1 in Section 2.3 of \cite{Homcables}, there exists a vertically simplified basis for $\CFKi{(11n_{57})}$. This means that we may assume that the vertical differential $\partial_{\textrm{vert}}$ cancels the basis elements in pairs, except for the one basis element that generates the vertical homology. Since $b,c$, and $d$ are the only basis elements that have homological grading 0, one of them must generate the vertical homology. Without loss of generality let $c$ be this element. Then, either $\partial_{\textrm{vert}}(a) = b$ or $\partial_{\textrm{vert}}(a) = d$. We choose the first option. Similarly, we choose \[\partial_{\textrm{vert}}(j) = k, \partial_{\textrm{vert}}(l) = n, \partial_{\textrm{vert}}(m) =o, \:\textrm{and}\: \partial_{\textrm{vert}}(p) = q.\] This results in the incomplete complex shown in Figure \ref{fig:cfk57_2}.

Furthermore, we note that by restrictions from the homological grading, we must have 
\begin{equation*}
  \partial(a) = b, \partial(b) = \epsilon_1 Ua \quad\textrm{and}\quad  \partial(c) = \epsilon_2 Ua, 
\end{equation*}
\noindent where $\epsilon_1, \epsilon_2 \in \{0,1\}.$ Then, 
\begin{equation*}
    \partial^2 (b) = \partial (\epsilon_1 Ua) = \epsilon_1 Ub = 0 \quad\textrm{and}\quad \partial^2 (c) = \partial (\epsilon_2 Ua) = \epsilon_2 Ub = 0,
\end{equation*}

\noindent which imply that $\epsilon_1 = \epsilon_2 = 0.$

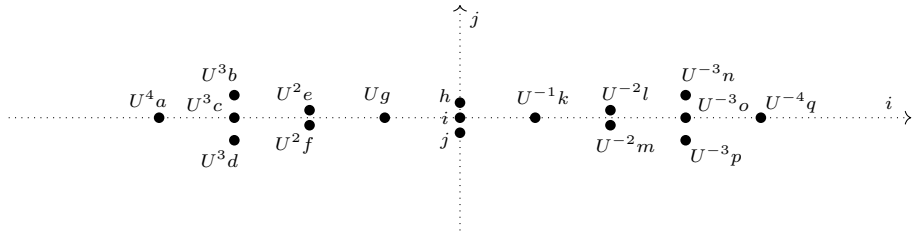
\begin{figure}
    \centering
    \begin{tikzpicture}\tikzstyle{every node}=[font=\tiny] 
    \path[->][dotted](0,-1.5)edge(0,1.5);
    \path[->][dotted](-6,0)edge(6,0);
    \node() at (.2,1.3){$j$};
    \node() at (5.7,.2){$i$};
    
	\fill(-4,0)circle [radius=2pt];
    \node(8) at (-4.15,0.25){$U^4a$};     
    
    \fill(-3,0.3)circle[radius=2pt];
    \node(7) at (-3.2,0.6){$U^3b$};
    \fill(-3,0)circle [radius=2pt];
    \node(6) at (-3.4,0.2){$U^3c$};
    \fill(-3,-0.3)circle [radius=2pt];
    \node(5) at (-3.2,-0.55){$U^3d$};    
    
	\fill(-2,0.1)circle[radius=2pt];
    \node(4) at (-2.2,0.35){$U^2e$};
    \fill(-2,-0.1)circle [radius=2pt];
    \node(3) at (-2.2,-0.35){$U^2f$};    
    
	\fill(-1,0)circle [radius=2pt];
    \node(2) at (-1.1,0.3){$Ug$};    
    
    \fill(0,0.2)circle[radius=2pt];
    \node(-1) at (-0.2,0.3){$h$};
    \fill(0,0)circle [radius=2pt];
    \node(0) at (-0.2,0){$i$};
    \fill(0,-0.2)circle [radius=2pt];
    \node(1) at (-0.2,-0.3){$j$};
    
    \fill(1,0)circle [radius=2pt];
    \node(-2) at (1.1,0.3){$U^{-1}k$}; 
    
    \fill(2,0.1)circle[radius=2pt];
    \node(-4) at (2.2,0.35){$U^{-2}l$};
    \fill(2,-0.1)circle [radius=2pt];
    \node(-3) at (2.2,-0.35){$U^{-2}m$};
    
    \fill(3,0.3)circle[radius=2pt];
    \node(-7) at (3.3,0.6){$U^{-3}n$};
    \fill(3,0)circle [radius=2pt];
    \node(-6) at (3.5,0.2){$U^{-3}o$};
    \fill(3,-0.3)circle [radius=2pt];
    \node(-5) at (3.4,-0.5){$U^{-3}p$}; 
    
    \fill(4,0)circle [radius=2pt];
    \node(-8) at (4.4,0.2){$U^{-4}q$}; 

    \end{tikzpicture}
    \caption{The arrangement of the generators of $C\{j=0\}.$ Note that $U^{-3}n, U^{-3}o$, and $U^{-3}p$ all have homological grading 0.}
    \label{fig:horGen}
\end{figure}

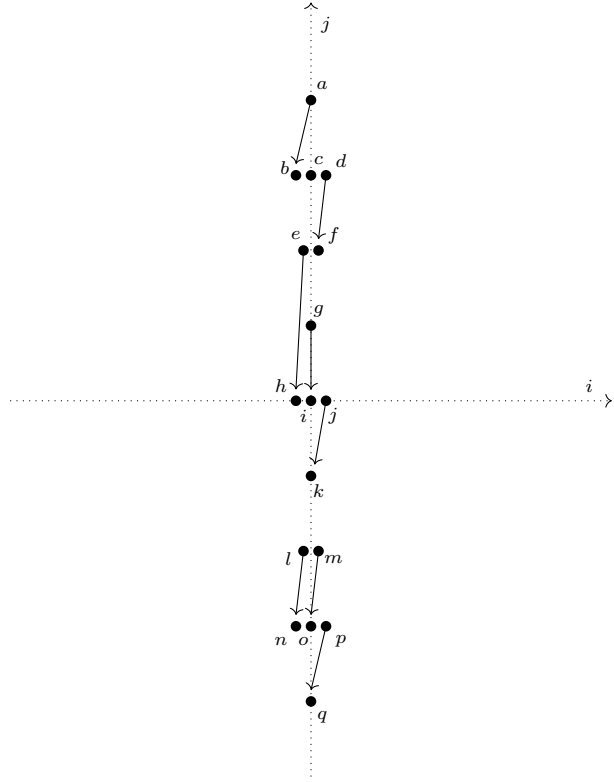
\begin{figure}
    \centering
    \begin{tikzpicture}\tikzstyle{every node}=[font=\tiny] 
    \path[->][dotted](0,-5)edge(0,5.3);
    \path[->][dotted](-4,0)edge(4,0);
    \node() at (.2,5){$j$};
    \node() at (3.7,.2){$i$};
    
	\fill(0,4)circle [radius=2pt];
    \node(8) at (0.15,4.2){$a$};
    \path[->](0,4)edge(-0.2,3.15);  
    
    \fill(-0.2,3)circle[radius=2pt];
    \node(7) at (-0.35,3.1){$b$};
    \fill(0,3)circle [radius=2pt];
    \node(6) at (0.1,3.2){$c$};
    \fill(0.2,3)circle [radius=2pt];
    \node(5) at (0.4,3.2){$d$};
    \path[->](0.2,3)edge(0.1,2.15);  
    
	\fill(-0.1,2)circle[radius=2pt];
    \node(4) at (-0.2,2.2){$e$};
    \path[->](-0.1,2)edge(-0.2,0.15);
    \fill(0.1,2)circle [radius=2pt];
    \node(3) at (0.3,2.2){$f$};    
    
	\fill(0,1)circle [radius=2pt];
    \node(2) at (0.1,1.2){$g$};
    \path[->](0,1)edge(0,0.15);  
    
    \fill(-0.2,0)circle[radius=2pt];
    \node(-1) at (-0.4,0.2){$h$};
    \fill(0,0)circle [radius=2pt];
    \node(0) at (-0.1,-0.2){$i$};
    \fill(0.2,0)circle [radius=2pt];
    \node(1) at (0.3,-0.2){$j$};
    \path[->](0.2,0)edge(0.05,-0.85);
    
    \fill(0,-1)circle [radius=2pt];
    \node(-2) at (0.1,-1.2){$k$}; 
    
    \fill(-0.1,-2)circle[radius=2pt];
    \node(-4) at (-0.3,-2.1){$l$};
    \path[->](-0.1,-2)edge(-0.2,-2.85);
    \fill(0.1,-2)circle [radius=2pt];
    \node(-3) at (0.3,-2.1){$m$};
    \path[->](0.1,-2)edge(0,-2.85);
    
    \fill(-0.2,-3)circle[radius=2pt];
    \node(-7) at (-0.4,-3.2){$n$};
    \fill(0,-3)circle [radius=2pt];
    \node(-6) at (-0.1,-3.2){$o$};
    \fill(0.2,-3)circle [radius=2pt];
    \node(-5) at (0.4,-3.2){$p$};
    \path[->](0.2,-3)edge(0,-3.85);
    
    \fill(0,-4)circle [radius=2pt];
    \node(-8) at (0.15,-4.2){$q$}; 
    \end{tikzpicture}
    \caption{An updated, and now complete picture of $C\{i=0\}$. All generators are paired up by the vertical differential, except for $c$, which generates the homology of the complex.}
    \label{fig:cfk57_3}
\end{figure}

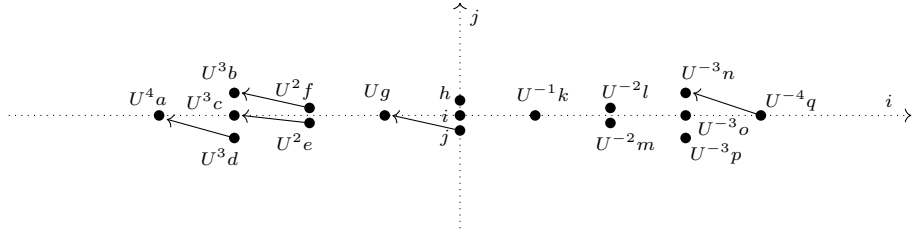
\begin{figure}
    \centering
    \begin{tikzpicture}\tikzstyle{every node}=[font=\tiny] 
    \path[->][dotted](0,-1.5)edge(0,1.5);
    \path[->][dotted](-6,0)edge(6,0);
    \node() at (.2,1.3){$j$};
    \node() at (5.7,.2){$i$};
    
	\fill(-4,0)circle [radius=2pt];
    \node(8) at (-4.15,0.25){$U^4a$};     
    
    \fill(-3,0.3)circle[radius=2pt];
    \node(7) at (-3.2,0.6){$U^3b$};
    \fill(-3,0)circle [radius=2pt];
    \node(6) at (-3.4,0.2){$U^3c$};
    \fill(-3,-0.3)circle [radius=2pt];
    \node(5) at (-3.2,-0.55){$U^3d$};
    \path[->](-3,-0.3)edge(-3.9,-0.05);    
    
	\fill(-2,0.1)circle[radius=2pt];
    \node(4) at (-2.2,0.35){$U^2f$};
    \path[->](-2,0.1)edge(-2.9,0.3);
    \fill(-2,-0.1)circle [radius=2pt];
    \node(3) at (-2.2,-0.35){$U^2e$};
    \path[->](-2,-0.1)edge(-2.9,-0);  
    
	\fill(-1,0)circle [radius=2pt];
    \node(2) at (-1.1,0.3){$Ug$};    
    
    \fill(0,0.2)circle[radius=2pt];
    \node(-1) at (-0.2,0.3){$h$};
    \fill(0,0)circle [radius=2pt];
    \node(0) at (-0.2,0){$i$};
    \fill(0,-0.2)circle [radius=2pt];
    \node(1) at (-0.2,-0.3){$j$};
    \path[->](0,-0.2)edge(-0.9,0);
    
    \fill(1,0)circle [radius=2pt];
    \node(-2) at (1.1,0.3){$U^{-1}k$}; 
    
    \fill(2,0.1)circle[radius=2pt];
    \node(-4) at (2.2,0.35){$U^{-2}l$};
    \fill(2,-0.1)circle [radius=2pt];
    \node(-3) at (2.2,-0.35){$U^{-2}m$};
    
    \fill(3,0.3)circle[radius=2pt];
    \node(-7) at (3.3,0.6){$U^{-3}n$};
    \fill(3,0)circle [radius=2pt];
    \node(-6) at (3.5,-0.15){$U^{-3}o$};
    \fill(3,-0.3)circle [radius=2pt];
    \node(-5) at (3.4,-0.5){$U^{-3}p$}; 
    
    \fill(4,0)circle [radius=2pt];
    \node(-8) at (4.4,0.2){$U^{-4}q$}; 
    \path[->](4,0)edge(3.1,0.3);

    \end{tikzpicture}
    \caption{An updated, but incomplete picture for $C\{j=0\}$ with the horizontal differential of $U^3d,U^2e,U^2f,j,$ and $U^{-4}q$ denoted by the arrows.} 
    \label{fig:cfkHorz_1}
\end{figure}

We now consider the horizontal differential $\partial_{\textrm{horz}}.$ The arrangement of the generators of $C\{j=0\}$ is shown in Figure \ref{fig:horGen}. We know that the homology of $C\{j=0\}$ is one-dimensional, and is generated by a linear combination of $U^{-3}n, U^{-3}o,$ or $U^{-3}p.$ By the restrictions on homological grading we must have $\partial_{\textrm{horz}}(g) = 0.$ Also, $U^4 a$ must be in the image of the horizontal differential, which implies that $\partial_{\textrm{horz}}(d) = Ua.$ Then, $$\partial^2 (d) = \partial(Ua) + \partial (\partial_{\textrm{vert}} d) = Ub + \partial (\partial_{\textrm{vert}} d),$$ implying that $\partial_{\textrm{vert}}(d) = f.$ Now we can choose $\partial_{\textrm{vert}}(e) = h$ and $\partial_{\textrm{vert}}(g) = i.$ We can also choose $\partial_{\textrm{horz}}(e) = Uc$ and $\partial_{\textrm{horz}}(f) = Ub$. This is because $U^3b$ and $U^3c$ must be in the image of the horizontal differential, so their preimages under this map must be linearly independent combinations of $U^2e$ and $U^2f$. Up to a change of basis the above choice is the only possibility. Since $Ug$ must be in the image of $\partial_{\textrm{horz}}$, we choose $\partial_{\textrm{horz}}(j) = Ug$. Since $U^{-4}q$ does not generate the homology of $C\{j=0\}$, we choose $\partial_{\textrm{horz}}(q) = Un.$ Given these calculations, the updated pictures for $C\{i=0\}$ and $C\{j=0\}$ are in Figure \ref{fig:cfk57_3} and Figure \ref{fig:cfkHorz_1}, respectively.

We now find the rest of the horizontal differential. First, we have $$\partial^2(e) = \partial(Uc + h) = \partial(h) = 0.$$ By restrictions on homological gradings, 

\begin{equation*}
    \partial (i) = \epsilon_1 Ug \quad\textrm{and}\quad \partial(j) = \epsilon_2 Ug + k,
\end{equation*}
\noindent where $\epsilon_1, \epsilon_2 \in \{0,1\}.$ Then, $$\partial^2 (i) = \partial (\epsilon_1 Ui) = 0,$$ which implies that $\epsilon_1 = 0.$ However, since $Ug$ is in the image of the horizontal, it must be that $\epsilon_2 = 1.$ Furthermore, $$\partial^2(j) = \partial (Ug + k) = Ui + \partial(k) = 0,$$ which implies that $\partial(k) = Ui.$ 
The images of $l,m,n,o,p$, and $q$ under $\partial$ are as follows.
\begin{align*}
    \partial (l) &= n + \kappa_1 U^2h + \kappa_2 U^2i + \kappa_3 U^2 j & \partial(n) &= \gamma_1 Ul + \gamma_2 Um\\
    \partial (m) &= o + \lambda_1 U^2h + \lambda_2 U^2i + \lambda_3 U^2 j & \partial (o) &= \beta_1 Ul + \beta_2 Um\\
    \partial (q) &= \epsilon_1 Un + \epsilon_2 Uo + \epsilon_3 Up & \partial (p) &= q + \alpha_1 Ul + \alpha_2 Um,
\end{align*} where all coefficients are either 0 or 1. Consider

\[\partial^2(l) = \partial(n + \kappa_1 U^2h + \kappa_2 U^2i + \kappa_3 U^2 j) = \gamma_1 Ul + \gamma_2 Um + \kappa_3(U^3g + Uk) = 0.\] Therefore, $\gamma_1 = \gamma_2 = \kappa_3 = 0.$ Similarly, \[\partial^2(m) = \partial(o + \lambda_1 U^2h + \lambda_2 U^2i + \lambda_3 U^2 j) = \beta_1 Ul + \beta_2 Um + \lambda_3(U^3g + Uk) = 0,\] which implies that $\beta_1 = \beta_2 = \lambda_3 = 0.$ So, both $\partial(n)$ and $\partial(o)$ are 0. Thus, \[\partial^2(q) = \partial(\epsilon_1 Un + \epsilon_2 Uo + \epsilon_3 Up) = \epsilon_3 (Uq + \alpha_1 U^2l + \alpha_2 U^2m),\] which implies that $\epsilon_3 = 0.$ Now, since $\partial_{\textrm{horz}}(q)$ must be nonzero we have that $\partial(q) = Un$, which is the only option up to relabeling and changing the basis. Lastly we consider \[\partial^2(p) = \partial (q + \alpha_1 Ul + \alpha_2 Um) = Un + \alpha_1 (Un + \kappa_1 U^3h + \kappa_2 U^3i) + \alpha_2 (Uo + \lambda_1 U^3h + \lambda_2 U^3i) = 0.\] It follows that $\alpha_2 = \kappa_1 = \kappa_2 = 0$ and $\alpha_1 = 1$. Additionally, $\lambda_1 = 1$, because $h$ does not generate the homology of $C\{j=0\}.$ Up to a change of basis we also have $\partial(m) = o + U^2h.$ We have determined the image of all generators of $\CFKi(11n_{57})$ under the boundary map. The resulting complex is shown in Figure \ref{fig:cfk57_final}.

\subsubsection{Finding the map $\iota_K$ for $11n_{57}$}
\label{sub:iotaK}

The complex associated to $11n_{57}$ is shown in Figure \ref{fig:cfk57_final}. All coefficients in the following discussion are either 0 or 1. Since $\iota_K$ is skew-filtered and grading-preserving, it must be of the form below. 
\begingroup
\allowdisplaybreaks
\begin{align*}
    &a \mapsto U^{-4}q 
    &&b \mapsto \beta_0U^{-3}p + \beta_1U^{-3}o + \beta_2U^{-3}n\\
    &c \mapsto \delta_0U^{-3}p + \delta_1U^{-3}o + \delta_2U^{-3}n
    &&d \mapsto \alpha_0U^{-3}p + \alpha_1U^{-3}o + \alpha_2U^{-3}n\\
    &e \mapsto \epsilon_0U^{-1}k + \epsilon_1U^{-2}m + \epsilon_2U^{-2}l
    &&f \mapsto \gamma_0U^{-1}k + \gamma_1U^{-2}m + \gamma_2U^{-2}l\\
    &g \mapsto U^{-1}k
    &&h \mapsto \eta_0j + \eta_1Ui + \eta_2h\\
    &i \mapsto \psi_0h + \psi_1i + \psi_2j
    && j \mapsto \zeta_0j + \zeta_1Ui + \zeta_2h \\
    & k \mapsto Ug
    &&l \mapsto \kappa_0U^2f + \kappa_1U^2e + \kappa_2U^2g\\
    &m \mapsto \theta_0U^2f + \theta_1U^2e + \theta_2U^2g
    &&n \mapsto \nu_0U^3d + \nu_1U^3b + \nu_2U^3c\\
    & o \mapsto \mu_0U^3d + \mu_1U^3b + \mu_2U^3c
    && p \mapsto \lambda_0U^3d + \lambda_1U^3b + \lambda_2U^3c\\
    &q \mapsto U^4a.
\end{align*}
\endgroup

Furthermore, since $\iota_K$ is a chain map, we can narrow down the options further to conclude that it has the form of the following map.

\begingroup
\allowdisplaybreaks
\begin{align*}
    &a \mapsto U^{-4}q
    &&b \mapsto U^{-3}n\\
    &c \mapsto \epsilon_1U^{-3}o + \epsilon_2U^{-3}n
    &&d \mapsto U^{-3}p + \alpha_1U^{-3}o + \alpha_2U^{-3}n\\
    &e \mapsto \epsilon_0U^{-1}k + \epsilon_1U^{-2}m + \epsilon_2U^{-2}l
    &&f \mapsto U^{-2}l \\
    &g \mapsto U^{-1}k
    &&h \mapsto \epsilon_0Ui + \epsilon_1h\\
    &i \mapsto i
    &&j \mapsto j + \zeta_1Ui + \zeta_2h\\
    &k \mapsto Ug
    && l \mapsto U^2f\\
    &m \mapsto \theta_0U^2f + \theta_1U^2e + \theta_2U^2g
    &&n \mapsto U^3b\\
    & o \mapsto \theta_0U^3b + \theta_1U^3c
    && p \mapsto U^3d + \lambda_1U^3b + \lambda_2U^3c\\
    &q \mapsto U^4a.
\end{align*}
\endgroup

Now we consider the requirement that $\iota_K^2 = \sigma$ up to chain homotopy equivalence, where $\sigma$ is the Sarkar involution. Since $\iota_K$ and $\sigma$ are both filtered, and because the boundary map $\partial$ reduces the homological grading by one, any map $H: CFK^{\infty}(11n_{57}) \rightarrow CFK^{\infty}(11n_{57})$ such that $\partial H + H\partial = \iota_K^2 + \sigma$ must be filtered and increase homological grading by one where it is nonzero. There is no such map except for the trivial map $H \equiv 0.$ Thus, we have $\iota_K^2 = \sigma.$ In this case we see that $\sigma$ is the identity map except at $d$, $j$, and $p$. At these points we have $\sigma(d) = d + b, \sigma(j) = j + i$, and $\sigma(p) = p + n.$ Applying this restriction results in a further simplification of the options for $\iota_K$, which is shown below.
\begingroup
\allowdisplaybreaks
\begin{align*}
    &a \mapsto U^{-4}q
    &&b \mapsto U^{-3}n\\
    &c \mapsto U^{-3}o + \epsilon_2U^{-3}n
    &&d \mapsto U^{-3}p + \alpha_1U^{-3}o + \alpha_2U^{-3}n\\
    &e \mapsto U^{-1}k + U^{-2}m + \epsilon_2U^{-2}l
    &&f \mapsto U^{-2}l\\
    &g \mapsto U^{-1}k
    &&h \mapsto Ui + h\\
    &i \mapsto i
    &&j \mapsto j + \zeta_1Ui + h\\
    &k \mapsto Ug
    &&l \mapsto U^2f \\
    &m \mapsto \epsilon_2U^2f + U^2e + U^2g
    &&n \mapsto U^3b\\
    &o \mapsto \epsilon_2U^3b + U^3c
    &&p \mapsto U^3d + \lambda_1U^3b + \alpha_1U^3c\\
    &q \mapsto U^4a.
\end{align*}
\endgroup

To simplify the possibilities for $\iota_K$ we make the following change of basis (and then immediately drop the primes):
    \begin{align*}
    &m' = m + \epsilon_2l \\
    &o' = o + \epsilon_2n \\
    &p' = p + \alpha_1o+\alpha_2n \\
    &j' = j + \zeta_1Ui.
    \end{align*}
As a result $\iota_K$ is of the form

\begingroup
\allowdisplaybreaks
\begin{align*}
    &a \mapsto U^{-4}q
    &&b \mapsto U^{-3}n\\
    &c \mapsto U^{-3}o
    &&d \mapsto U^{-3}p\\
    &e \mapsto U^{-1}k + U^{-2}m
    &&f \mapsto U^{-2}l\\
    &g \mapsto U^{-1}k
    &&h \mapsto Ui + h\\
    &i \mapsto i
    &&j \mapsto j + h\\
    &k \mapsto Ug
    &&l \mapsto U^2f \\
    &m \mapsto U^2e + U^2g
    &&n \mapsto U^3b\\
    &o \mapsto U^3c
    &&p \mapsto U^3d + U^3b\\
    &q \mapsto U^4a.
\end{align*}
\endgroup

All other maps that satisfy the necessary properties of $\iota_K$ are equivalent to this solution up to a change of basis. We also see that the boxes in the second and fourth quadrants interact only with each other, and that no other part of the complex interacts with them. Thus, they need not be considered in the calculation of the concordance invariants, since they can be removed as an equivariant summand.


\subsubsection{Finding $\underline{V}_0(11n_{57})$ and $\overline{V}_0(11n_{57})$}\label{sub:11n57invariants}

Using the map $\iota_K$, we may now compute $V_0(11n_{57})$ and the involutive concordance invariants $\overline{V_0}(11n_{57})$ and $\underline{V}_0(11n_{57})$.

\begin{figure}
    \centering
    \begin{tikzpicture}\tikzstyle{every node}=[font=\tiny]
    \fill(-1,2)circle [radius=2pt];
    \node at (-1.2,2.3){$U^3c$};
    \fill(0,2)circle [radius=2pt];
    \node at (0.2,2.3){$U^2e$};
    \path[->](0,2)edge(-0.9,2);
    \path[->](0,2)edge(0,0.1);
    
    \fill(0,0)circle[radius=2pt];
    \node at (-0.3,-0.2){$U^2h$};
    
    \fill(1,1)circle[radius=2pt];
    \node at (0.7,0.8){$Uh$};
    
    \fill(2,2)circle[radius=2pt];
    \node at (1.8,1.8){$h$};

    \fill(2,0)circle [radius=2pt];
    \node at (2.2,0.2){$m$};
    \path[->](2,0)edge(0.1,0);
    \path[->](2,0)edge(2,-0.9);
    \fill(2,-1)circle [radius=2pt];
    \node at (2.2,-0.8){$o$};
    \end{tikzpicture}
    \caption{$A_0^-$ for the staircase subcomplex of $\CFKi{(11n_{57})}$ consists of all non-negative $U$-translates of the staircase shown, along with the elements $h$ and $Uh$.}
    \label{fig:11n_57_A0}
\end{figure}
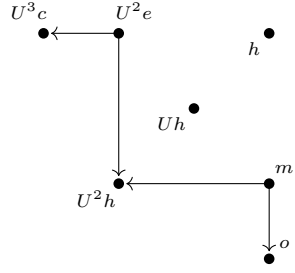

To find $V_0$ it is enough to consider the staircase in $\CFKi{(11n_{57})}$, since the square in the first quadrant has trivial homology. The complex $A_0^-$ for the staircase is shown in Figure \ref{fig:11n_57_A0}. The homology of $A_0^-$ is $\F_2[U] \langle [h] \rangle$. Thus, $V_0(11n_{57})$ is $-\frac{1}{2}$ times the homological grading of $h$, which yields $$V_0(11n_{57}) = -\frac{1}{2}(-2) = 1.$$

\begin{figure}
    \centering
    \begin{tikzpicture}\tikzstyle{every node}=[font=\tiny]
    \fill(-1,2)circle [radius=2pt];
    \node at (-1.2,2.3){$U^3c$};
    \fill(0,2)circle [radius=2pt];
    \node at (0.2,2.3){$U^2e$};
    \path[->](0,2)edge(-0.9,2);
    \path[->](0,2)edge(0,0.1);
    
    \fill(0,0)circle[radius=2pt];
    \node at (-0.3,-0.2){$U^2h$};
    
    \fill(0.98,1.15)circle[radius=2pt];
    \node at (0.85,1.4){$Uh$};

    \fill(2,0)circle [radius=2pt];
    \node at (2.2,0.2){$m$};
    \path[->](2,0)edge(0.1,0);
    \path[->](2,0)edge(2,-0.9);
    \fill(2,-1)circle [radius=2pt];
    \node at (2.2,-0.8){$o$};
    
    \fill(0.15,0.95)circle [radius=2pt];
    \node at (0.3,1.2){$U^2g$};
    \path[->](0.15,0.95)edge(0.15,0.25);  
    \fill(0.15,0.15)circle [radius=2pt];
    \node at (0.55,0.4){$U^2i$};
    \fill(0.95,0.95)circle [radius=2pt];
    \node at (1.25,0.8){$Uj$};
    \path[->](0.95,0.95)edge(0.95,0.25);
    \path[->](0.95,0.95)edge(0.25,0.95);
    \fill(0.95,0.15)circle [radius=2pt];
    \node at (1.25,0.3){$Uk$};
    \path[->](0.95,0.15)edge(0.25,0.15);
    
	\fill(1.15,1.9)circle [radius=2pt];
	\node at (1.25,2.2){$Ug$};
    \path[->](1.15,1.9)edge(1.15,1.2);  
    \fill(1.15,1.1)circle [radius=2pt];
    \node at (1.5,1.3){$Ui$};
    \fill(1.95,1.9)circle [radius=2pt];
    \node at (2.15,1.7){$j$};
    \path[->](1.95,1.9)edge(1.95,1.2);
    \path[->](1.95,1.9)edge(1.25,1.9);
    \fill(1.95,1.1)circle [radius=2pt];
    \node at (2.15,1.3){$k$};
    \path[->](1.95,1.1)edge(1.25,1.1);
    
    \fill(1.98,2.1)circle[radius=2pt];
    \node at (1.85,2.3){$h$};
    \fill(2.15,2.05)circle[radius=2pt];
    \node at (2.3,2.2){$i$};
    \end{tikzpicture}
    \vspace{0.7in}
    \caption{$A_0^-$ for the subcomplex of $\CFKi{(11n_{57})}$ made up of the staircase and box in the first quadrant consists of all non-negative $U$-translates of the staircase shown, all positive $U$-translates of the box, and the elements $i$,$h$, and $Uh$.}
    \label{fig:11n_57A0-big}
\end{figure}
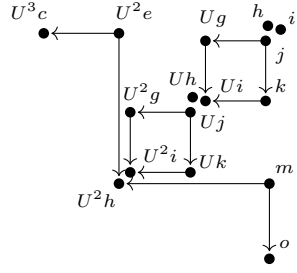

Now we calculate the involutive concordance invariants $\underline{V}_0$ and $\overline{V}_0$. To do so, we need to examine the homology of the mapping cone \[AI_0^- = \textrm{Cone}(A_0^- \xrightarrow{Q(\iota_{11n_{57}} + \textrm{Id})} QA_0^-[-1]).\] First, we consider $A_0^-$. The complex $A_0^-$ associated to the subcomplex of $\CFKi{(11n_{57})}$ made up of the staircase and box in quadrant 1 is shown in Figure \ref{fig:11n_57A0-big}. We also consider the corresponding complex $QA_0^-$.

To simplify the picture and the subsequent calculations we choose the following change of bases (and immediately drop the primes):
    \begin{align*}
    &c' = c + U^{-3}o
    &&Qc' = Qc + QU^{-3}o\\
    &g' = g + U^{-1}k
    &&Qg' = Qg + QU^{-1}k\\
    &e' = e + U^{-2}m
    &&Qe' = Qe + QU^{-2}m.
    \end{align*}

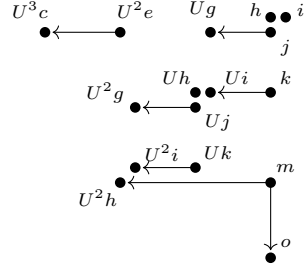
\begin{figure}
    \centering
    \begin{tikzpicture}\tikzstyle{every node}=[font=\tiny]
    \fill(-1,2)circle [radius=2pt];
    \node at (-1.2,2.3){$U^3c$};
    \fill(0,2)circle [radius=2pt];
    \node at (0.2,2.3){$U^2e$};
    \path[->](0,2)edge(-0.9,2);
    
    \fill(0,0)circle[radius=2pt];
    \node at (-0.3,-0.2){$U^2h$};
    
    \fill(2,2.2)circle[radius=2pt];
    \node at (1.8,2.3){$h$};
    \fill(2.2,2.2)circle[radius=2pt];
    \node at (2.4,2.3){$i$};

    \fill(2,0)circle [radius=2pt];
    \node at (2.2,0.2){$m$};
    \path[->](2,0)edge(0.1,0);
    \path[->](2,0)edge(2,-0.9);
    \fill(2,-1)circle [radius=2pt];
    \node at (2.2,-0.8){$o$};
    
    \fill(0.2,1)circle [radius=2pt];
    \node at (-0.2,1.2){$U^2g$};
    \fill(0.2,0.2)circle [radius=2pt];
    \node at (0.55,0.4){$U^2i$};
    \fill(1,1)circle [radius=2pt];
    \node at (1.3,0.8){$Uj$};
    \path[->](1,1)edge(0.3,1);
    \fill(1,0.2)circle [radius=2pt];
    \node at (1.3,0.4){$Uk$};
    \path[->](1,0.2)edge(0.3,0.2);
    
	\fill(1,1.2)circle[radius=2pt];
    \node at (0.75,1.4){$Uh$};    
    
	\fill(1.2,2)circle [radius=2pt];
	\node at (1.1,2.3){$Ug$};
    \fill(1.2,1.2)circle [radius=2pt];
    \node at (1.55,1.4){$Ui$};
    \fill(2,2)circle [radius=2pt];
    \node at (2.2,1.8){$j$};
    \path[->](2,2)edge(1.3,2);
    \fill(2,1.2)circle [radius=2pt];
    \node at (2.2,1.4){$k$};
    \path[->](2,1.2)edge(1.3,1.2);
    \end{tikzpicture}   \vspace{0.7in}
    \caption{$A_0^-$ for the subcomplex of $\CFKi{(11n_{57})}$ made up of the staircase and box in the first quadrant after a change of basis.}
    \label{fig:11n_57A0-basis}
\end{figure}

The resulting complex is shown in Figure \ref{fig:11n_57A0-basis}. Now, we look at the map $Q(\iota_{11n_{57}}+\textrm{Id})$, where $\iota_{11n_{57}}$ is given by the map described at the end of Section \ref{sub:iotaK}. After the change of bases, the restriction of $Q(\iota_{11n_{57}} + \textrm{Id})$ to the staircase and box in quadrant 1 is given by:

\begingroup
\begin{align*}
    &c \mapsto 0
    &&e \mapsto Qg\\
    &g \mapsto 0
    &&h \mapsto Qi\\
    &i \mapsto 0
    &&j \mapsto Qh\\
    &k \mapsto QUg
    &&m \mapsto QU^2e + QU^2g + QUk\\
    &o \mapsto QU^3c.
\end{align*}
\endgroup

Putting all the information together, the complete map on $AI_0^-$ is given by:

\begingroup
\begin{align*}
    &U^3c \mapsto 0
    &&QU^3c \mapsto 0\\
    &U^2e \mapsto U^3c + QU^2g
    &&QU^2e \mapsto QU^3c\\
    &Ug \mapsto 0
    &&QUg \mapsto 0\\
    &h \mapsto Qi
    &&Qh \mapsto 0\\
    &i \mapsto 0
    &&Qi \mapsto 0\\
    &j \mapsto Ug + Qh
    &&Qj \mapsto QUg\\
    &k \mapsto Ui + QUg
    &&Qk \mapsto QUi\\
    &m \mapsto U^2h + o + QU^2e + QU^2g + QUk
    &&Qm \mapsto QU^2h + Qo\\
    &o \mapsto QU^3c
    &&Qo \mapsto 0.
\end{align*}
\endgroup

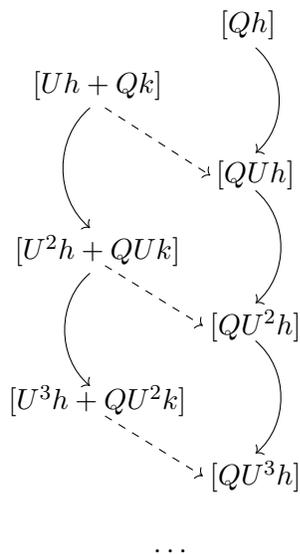
\begin{figure}
    \centering
    \begin{tikzpicture}
    \node(0)at(4,2){$[i]$};
    \node(1)at(0,0.2){$[Uh + Qk]$};
    \path[->][bend right = 50](-.1,-.1)edge(-.1,-1.7);
    \path[->][dashed](.1,-.1)edge(1.5,-1);
    \node(2)at(0,-2){$[U^2h + QUk]$};
    \path[->][bend right = 50](-.1,-2.3)edge(-.1,-3.8);
    \path[->][dashed](.1,-2.2)edge(1.4,-3);
    \node(3)at(0,-4){$[U^3h + QU^2k]$};
    \path[->][dashed](.1,-4.2)edge(1.4,-5);
    \node(4)at(2,1){$[Qh]$};
    \path[->][bend left = 50](2.1,0.7)edge(2.1,-0.7);
    \node(5)at(2.1,-1){$[QUh]$};
    \path[->][bend left = 50](2.1,-1.2)edge(2.1,-2.7);
    \node(6)at(2.1,-3){$[QU^2h]$};
    \path[->][bend left = 50](2.1,-3.2)edge(2.1,-4.7);
    \node(7)at(2.1,-5){$[QU^3h]$};
    \node(8)at(1,-6){$\bf \cdots$};
\end{tikzpicture}
    \caption{$H_*(AI_0^-)$ for the knot $11n_{57}$ can be described by two linked towers and the stand-alone subspace $[i]$.}
    \label{fig:11n57_towers}
\end{figure}

Now, we can extract the homology of $AI_0^-$, and form two towers as shown in Figure \ref{fig:11n57_towers}. We examine the homological gradings of the tops of the two towers. The homological grading of $[Qk]$ is $-3$, while the grading of $[Qh]$ is $-2$. Thus the involutive concordance invariants are $\underline{V}_0 = 2$ and $\overline{V}_0 = 1$.

\subsection{Finding the involutive concordance invariants for $10_{161}$.}\label{sub:10_161ex}
We now provide a second example, in which we extract $\CFKi(10_{161})$ from a (1,1) Heegaard diagram for the knot, then compute the automorphism $\iota_K$ and the involutive concordance invariants.
\subsubsection{Finding $\CFKi(10_{161})$ using a $(1,1)$ Heegaard diagram for $10_{161}$.}\label{sub:10_161cfk}

\begin{figure}
    \centering
    \includegraphics[width=.7\linewidth]{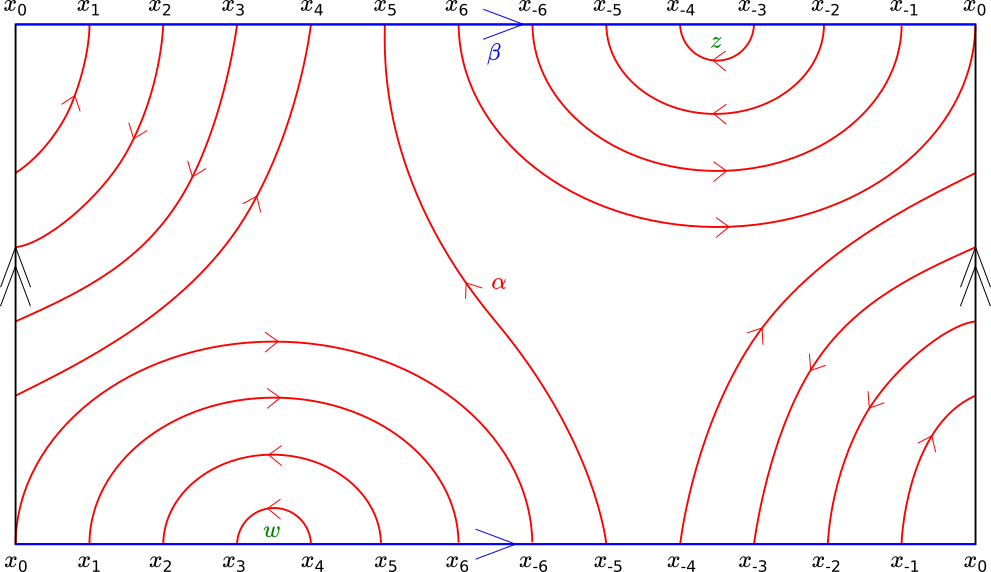}
    \caption{The (1,1)-diagram for $10_{161}$ has 13 intersection points, labelled by $x_i$ for $-6 \leq i \leq 6$, between the loops $\alpha$ and $\beta$.}
    \label{fig:10_161(1,1)}
\end{figure}
\begin{figure}
    \centering
    \includegraphics[width=.7\linewidth]{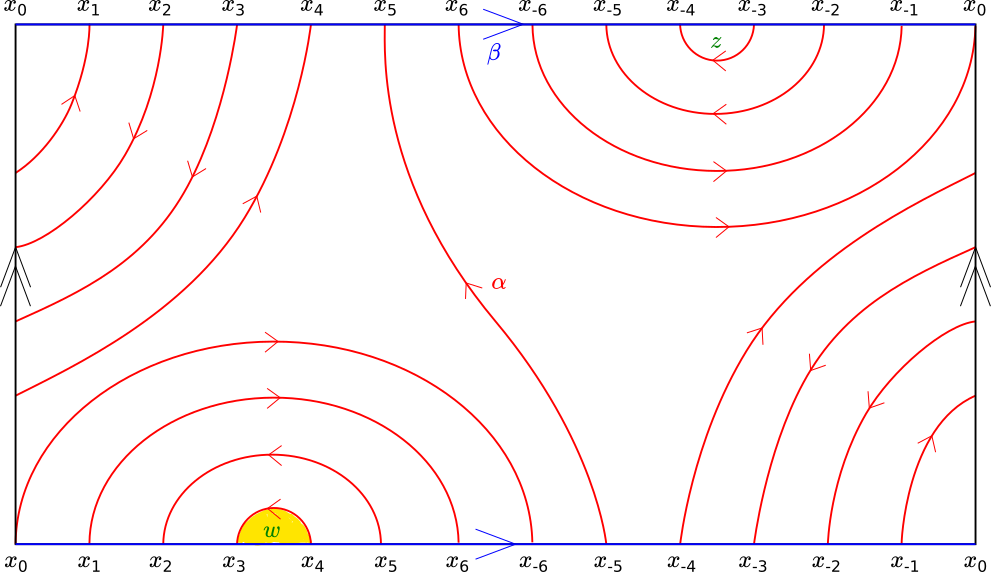}
    \caption{The disk highlighted in yellow is a  bigon from $x_4$ to $x_3$ containing one copy of $w$; the corresponding arrow in the chain complex runs from $x_4$ to $Ux_3$.}
    \label{fig:Bigon}
\end{figure}
\begin{figure}
    \centering
    \includegraphics[width=.7\linewidth]{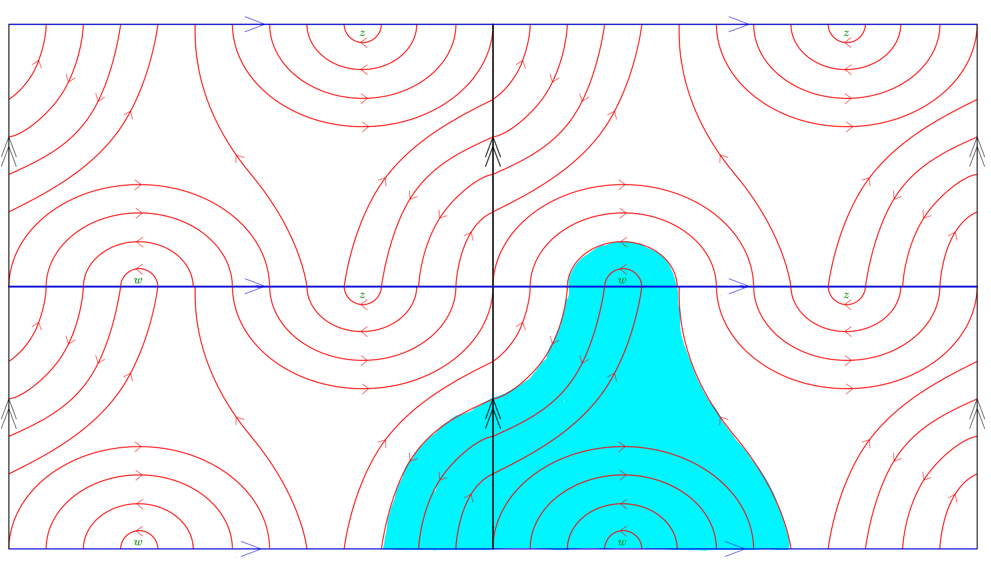}
    \caption{The disk highlighted in blue is a bigon from $x_{-5}$ to $x_{-3}$ containing two copies of $w$; the corresponding arrow in the chain complex runs from $x_{-5}$ to $U^2x_{-3}$.}
    \label{fig:Bigon2}
\end{figure}

A $(1,1)$ Heegaard diagram for $10_{161}$, computed in \cite[Section 3]{Racz}, is shown in Figure \ref{fig:10_161(1,1)}. We may compute $\CFKi{(10_{161})}$ from Figure \ref{fig:10_161(1,1)} by identifying bigons in this Heegaard diagram. Figure \ref{fig:Bigon} and Figure \ref{fig:Bigon2} show examples of bigons. Each bigon determines a component of the boundary map for $\CFKi{(10_{161})}$. For the bigon in Figure \ref{fig:Bigon}, the curves $\alpha$ and $\beta$ intersect at $x_3$ and $x_4$. When the bigon is oriented such that the segment of $\beta$ on its boundary is on the right, $x_4$ is above $x_3$. Furthermore, the bigon contains one copy of the point $w$. We conclude that a component of the boundary map for $\CFKi{(10_{161})}$ gives an arrow from $x_4$ to $Ux_3$ in the total differential, as shown in Figure \ref{fig:10_161cfkbig}.

\begin{figure}[ht] 
    \centering
    \begin{tikzpicture}\tikzstyle{every node}=[font=\tiny]
    \path[->][dotted](0,-4)edge(0,4.3);
    \path[->][dotted](-6,0)edge(6,0);
    \node(0) at (-4.6, 2.1){$U^2x_{-4}$};
    \node(1) at (-2.4,2.3){$Ux_1$};
    \node(2) at (0.3,2.1){$x_6$};
    \node(3) at (-4.6, -0.2){$U^2x_{-3}$};
    \node(4) at (-2.4, -0.2){$Ux_2$};
    \node(5) at (-0.6,-0.2){$x_5$};
    \node(6) at (0.3,0.5){$x_0$};
    \node(7) at (2.7, 0.3){$U^{-1}x_{-6}$};
    \node(8) at (0.4, -.5){$x_{-5}$};
    \node(9) at (0.4,-2.3){$x_{-2}$};
    \node(10) at (2.7,-1.9){$U^{-1}x_{-1}$};
    \node(11) at (0.3, -4.3){$x_3$};
    \node(12) at (2.6,-4.3){$U^{-1}x_4$};
    \fill(-4,2)circle [radius=2pt];
    \fill(-2,2)circle [radius=2pt];
    \fill(0,2)circle [radius=2pt];
    \fill(-4,0)circle [radius=2pt];
    \fill(-2,0)circle [radius=2pt];
    \fill(-0.2,0)circle [radius=2pt];
    \fill(0,0.2)circle [radius=2pt];
    \fill(2,0)circle [radius=2pt];
    \fill(0,-0.2)circle [radius=2pt];
    \fill(0,-2)circle [radius=2pt];
    \fill(2,-2)circle [radius=2pt];
    \fill(0,-4)circle [radius=2pt];
    \fill(2,-4)circle [radius=2pt];
    
    \path[->](-4,2)edge(-4,0.1);
    \path[->](-2,2)edge(-2,0.1);
    \path[->](-.2,0)edge(-1.9,0);
    \path[->](0,-.2)edge(0,-1.9);
    \path[->](2,-2)edge(0.1,-2);
    \path[->](2,0)edge(2,-1.9);
    \path[->](2,-4)edge(0.1,-4);
    \path[->](0,2)edge(-1.9,2);
    \path[->][bend right =7](2,0)edge(0.1,0.2);
    \path[->][bend left = 7](2,0)edge(2,-3.9);
    \path[->][bend left =7](2,0)edge(0.1,-.2);
    \path[->][bend left =7](0,2)edge(0,0.3);
    \path[->][bend right =7](0,2)edge(-3.9,2);
    \path[->][bend right =10](0,2)edge(-.2,0.1);
    \path[->][bend right =15](-0.05,0.1)edge(0,-3.9);
    \path[->][bend right =20](-.2,0)edge(-0.1,-3.95);
    \path[-][bend right =7](0,0.2)edge(-0.15,0.25);
    \path[-][bend right =5](-0.25,0.27)edge(-1.85,0.33);
    \path[->][bend right =5](-2.15,0.33)edge(-3.9,0.1);
    %
    \path[-][bend left = 15](0,-.2)edge(-0.1,-.25);
    \path[-][bend left = 15](-0.2,-.27)edge(-0.25,-.28);
    \path[->][bend left = 15](-0.35,-.3)edge(-3.9,-.07);
    \end{tikzpicture}
    \caption{$\CFKi{(10_{161})}$ prior to a change of basis.}
    \label{fig:10_161cfkbig}
\end{figure}
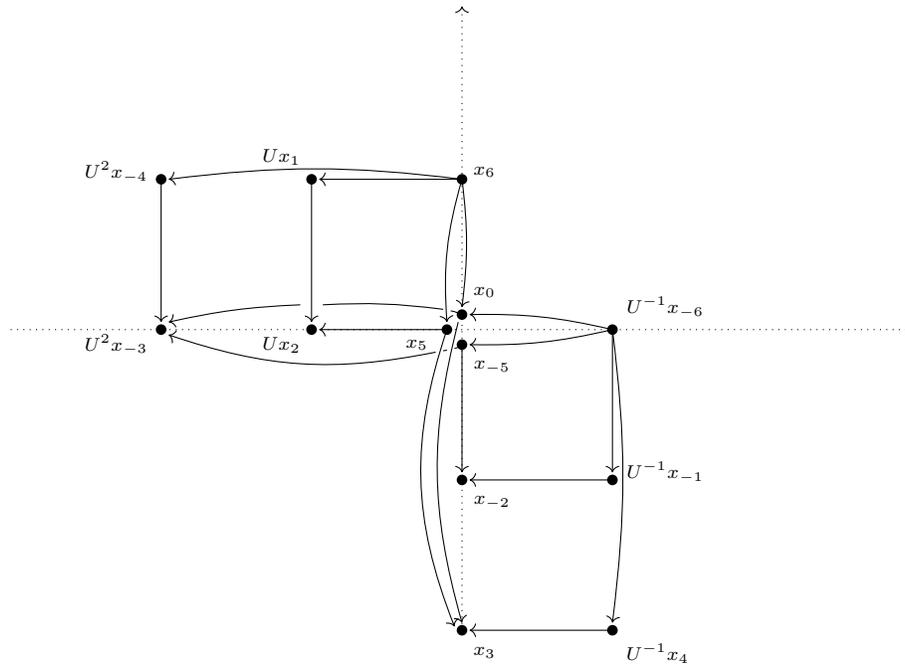

\begin{figure}
    \centering
    \begin{tikzpicture}\tikzstyle{every node}=[font=\tiny]
    \path[->][dotted](0,-4)edge(0,4.3);
    \path[->][dotted](-6,0)edge(6,0);
    \node(0) at (-4.6, 2.1){$U^2x_{-4}$};
    \node(1) at (-2.4,2.3){$Ux_1+U^2x_{-4}$};
    \node(2) at (0.3,2.1){$x_6$};
    \node(3) at (-4.6, -0.2){$U^2x_{-3}$};
    \node(4) at (-2.95, 0.5){$Ux_2+U^2x_{-3}$};
    \node(5) at (-0.4,-0.2){$x_0$};
    \node(6) at (0.6,0.5){$x_0+x_5$};
    \node(7) at (2.7, -0.3){$U^{-1}x_{-6}$};
    \node(8) at (0.9, -.3){$x_0+x_{-5}$};
    \node(9) at (0.9,-2.3){$x_{-2}+x_3$};
    \node(10) at (2.7,-2.3){$U^{-1}x_{-1}$};
    \node(11) at (0.3, -4.3){$x_3$};
    \node(12) at (2.6,-4.3){$U^{-1}x_4$};
    \fill(-4,2)circle [radius=2pt];
    \fill(-2,2)circle [radius=2pt];
    \fill(0,2)circle [radius=2pt];
    \fill(-4,0)circle [radius=2pt];
    \fill(-2,0.2)circle [radius=2pt];
    \fill(0,0)circle [radius=2pt];
    \fill(0,0.2)circle [radius=2pt];
    \fill(2,0)circle [radius=2pt];
    \fill(0.2,0)circle [radius=2pt];
    \fill(0.2,-2)circle [radius=2pt];
    \fill(2,-2)circle [radius=2pt];
    \fill(0,-4)circle [radius=2pt];
    \fill(2,-4)circle [radius=2pt];
    \path[->](-4,2)edge(-4,0.1);
    \path[->](-2,2)edge(-2,0.3);
    \path[->](0,0.2)edge(-1.9,0.2);
    \path[->](0.2,0)edge(0.2,-1.9);
    \path[->](2,-2)edge(0.3,-2);
    \path[->](2,0)edge(2,-1.9);
    \path[->](2,-4)edge(0.1,-4);
    \path[->](0,2)edge(-1.9,2);
    \path[->](2,0)edge(0.3,0);
    \path[->](0,2)edge(0,0.3);
    \path[->](0,0)edge(0,-3.9);
    \path[->](0,0)edge(-3.9,0);
    \end{tikzpicture}
    \caption{After the change of basis shown, $\CFKi{(10_{161})}$ consists of a staircase and two boxes.}
    \label{fig:10_161cfk}
\end{figure}
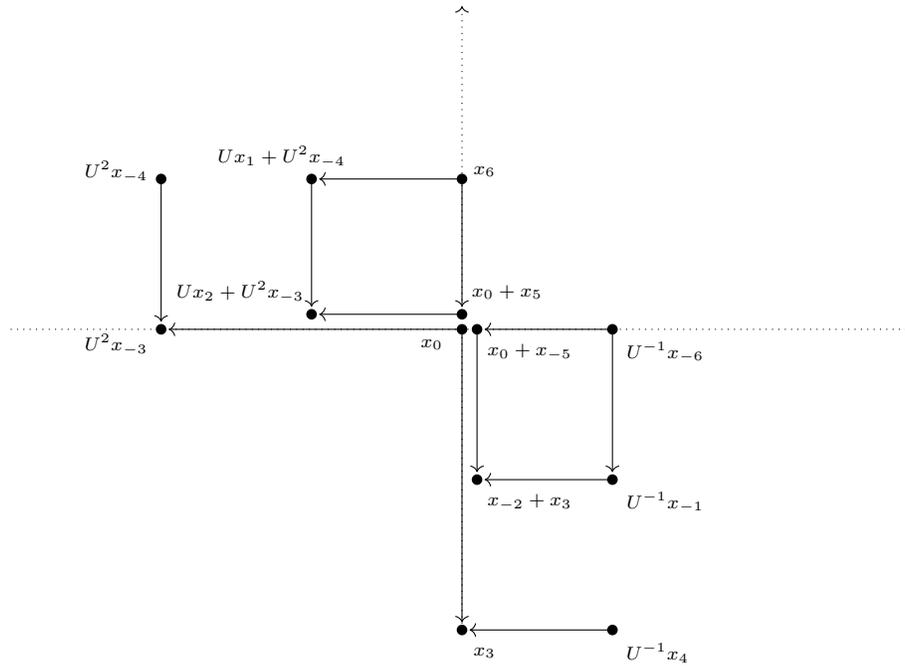

\subsubsection{Finding the map $\iota_K$ for $10_{161}$}\label{sub:10_161iota}

Now, we compute the automorphism $\iota_K$. We make the following change of basis, shown in Figure \ref{fig:10_161cfk}.
\begin{align*}
    x_1' &= x_1+Ux_{-4} \\
    x_2' &= x_2+Ux_{-3} \\
    x_5' &= x_5+x_0 \\
    x_{-1}' &= x_{-1}+x_4 \\
    x_{-2}' &= x_{-2}+x_3\\
    x_{-5}' &= x_{-5}+x_0
\end{align*}
%
To reduce the possibilities for $\iota_K$, we follow the same procedure as in Section \ref{sub:iotaK}, where we use that $\iota_K$ is skew-filtered, grading-preserving, a chain map, and squares to the Sarkar involution up to filtered chain homotopy. This process yields two possibilities for $\iota_K$. One is the map below, which we will denote $\iota_K$.
%

\begin{align*}
    &x_0 \mapsto x_0 && \hphantom{c} \\
    &x_1' \mapsto U^{-2}x_{-1}' && x_{-1}' \mapsto U^2x_1'\\
    &x_2' \mapsto U^{-1}x_{-2}' && x_{-2}' \mapsto Ux_2'\\
    &x_3 \mapsto U^2x_{-3} && x_{-3} \mapsto U^{-2}x_3\\
    &x_4 \mapsto U^3x_{-4}' &&  x_{-4} \mapsto U^{-3}x_4\\
    &x_5' \mapsto x_{-5}' && x_{-5}' \mapsto x_{5}'\\
    &x_6 \mapsto U^{-1}x_{-6} && x_{-6} \mapsto Ux_6 + x_2'
\end{align*}
%
The other is the map below, which we will denote $\iota_K'$.
\begin{align*}
    &x_0 \mapsto x_0 && \hphantom{c} \\
    &x_1' \mapsto U^{-2}x_{-1}' && x_{-1}' \mapsto U^2x_1'\\
    &x_2' \mapsto U^{-1}x_{-2}' && x_{-2}' \mapsto Ux_2'\\
    &x_3 \mapsto U^2x_{-3} && x_{-3} \mapsto U^{-2}x_3\\
    &x_4 \mapsto U^3x_{-4}' &&  x_{-4} \mapsto U^{-3}x_4\\
    &x_5' \mapsto x_{-5}' && x_{-5}' \mapsto x_{5}'\\
    &x_6 \mapsto U^{-1}x_{-6}+U^{-1}x_3 && x_{-6} \mapsto Ux_6 + x_2'
\end{align*}
\noindent The maps $\iota_K$ and $\iota_K'$ are equivalent up to a skew-filtered chain homotopy, the map $G$ below. 
\begin{align*}
    &x_1' \mapsto U^{-1}x_3 \\
    &x_i \mapsto 0 \text{ for all other generators } x_i 
\end{align*}
\noindent for which we have
\begin{align*}
    & \del G(x_1') + G \del (x_1')=0\\
    & \del G(x_6) + G \del (x_6)=G(x_1')=U^{-1}x_3.\end{align*}

%
%
%
%

\noindent This is exactly the map $\iota_K - \iota_K'$.

Thus, the two maps $\iota_K$ and $\iota_K'$ are related by a skew-filtered chain homotopy $G$, and we may use either in our computation of the involutive concordance invariants for $10_{161}$.

\subsubsection{Finding $\underline{V}_0(10_{161})$ and $\overline{V}_0(10_{161})$}\label{sub:10_161invariants}

Now, we compute the involutive concordance invariants $\overline{V}_0$ and $\underline{V}_0$.

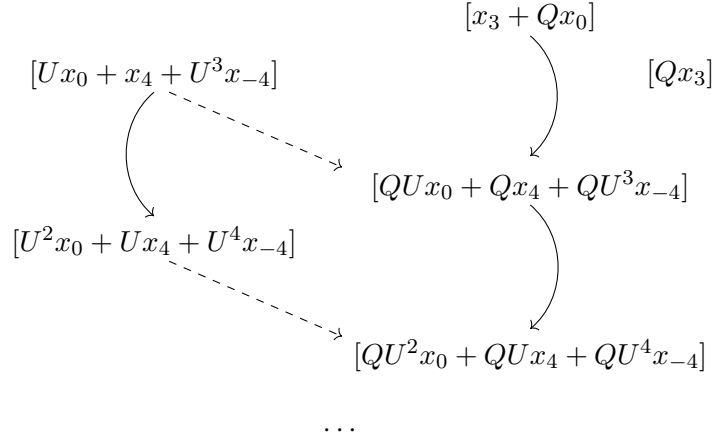
\begin{figure}
    \centering
    \begin{tikzpicture}
    \node(top)at(8,4.75){$[Qx_3]$};
    \node(1)at(1,4.75){$[Ux_0+x_4+U^3x_{-4}]$};
    \path[->][dashed](1.2,4.5)edge(3.5,3.5);
    \node(2)at(1,2.5){$[U^2x_0+Ux_4+U^4x_{-4}]$};
    \path[->][dashed](1.2,2.25)edge(3.5,1.25);
    \node(3)at(6,5.5){$[x_3+Qx_0]$};
    \node(4)at(6,3.25){$[QUx_0+Qx_4+QU^3x_{-4}]$};
    \node(5)at(6,1){$[QU^2x_0+QUx_4+QU^4x_{-4}]$};
    \node(6)at(3.5,0){$\bf \cdots$};
    
    \path[->][bend left = 50](6,5.25)edge(6,3.65);
    \path[->][bend left = 50](6,3)edge(6,1.35);
    \path[->][bend right = 50](1,4.5)edge(1,2.85);

\end{tikzpicture}
    \caption{$H_*(AI_0^-)$ for the knot $10_{161}$ can be written as two towers related by the action of $Q$ and the stand-alone subspace $[Qx_3]$.}
    \label{fig:10_161_towers}
\end{figure}
As in the computation of the involutive concordance invariants of $11n_{57}$ in Section \ref{sub:11n57ex}, we examine the homology of the mapping cone  \[AI_0^- = \textrm{Cone}(A_0^- \xrightarrow{Q(\iota_{10_{161}} + \textrm{Id})} QA_0^-[-1]).\]

We consider the complex $A_0^-$ associated to the subcomplex in the first quadrant and the corresponding complex $QA_0^-$. We extract the homology of $AI_0^-$ to form two towers, as shown in Figure \ref{fig:10_161_towers}, then examine the homological gradings of the tops of the two towers. We see that the homological grading of $[Qx_3]$ is $1$, the homological grading of $[x_3+Qx_0]$ is $2$, and the homological grading of $[Ux_0+x_4+U^3x_{-4}]$ is $1$. Thus, $\overline{V_0}(10_{161})=-1$ and
$\underline{V}_0(10_{161})=0$.

\FloatBarrier

\section{Results}\label{sec:results}
\subsection{Table of Invariants}\label{subsec:table} In Table \ref{tab:invol} we show the results of our computations for the involutive concordance invariants for all 10- and 11-crossing (1,1)-knots that are neither thin nor $L$-space.

\begin{table}[htp]
    \centering
    \begin{tabular}{|c|c|c|c|c|c|c|}
    \hline
    $\textbf{Knot } \bm{K}$ & $\bm{V_0(K)}$ & $\bm{\underline{V}_0(K)}$ & $\bm{\overline{V}_0(K)}$ & $\bm{V_0(\overline{K})}$ & $\bm{\underline{V}_0(\overline{K})}$ & $\bm{\overline{V}_0(\overline{K})}$ \\ \hline
    $10_{128}$    & 1          & 1                      & 1                     & 0                          & 0                                      & $-1$                                  \\ \hline
    $10_{132}$    & 0          & 0                      & $-1$                  & 1                          & 1                                      & 1                                     \\ \hline
    $10_{136}$    & 0          & 0                      & 0                     & 0                          & 0                                      & 0                                     \\ \hline
    $10_{139}$    & 2          & 2                      & 1                     & 0                          & 0                                      & $-2$                                  \\ \hline
    $10_{145}$    & 0          & 0                      & $-1$                  & 1                          & 1                                      & 1                                     \\ \hline
    $10_{161}$    & 0          & 0                      & $-1$                  & 1                          & 1                                      & 1                                     \\ \hline
    $11n_{12}$    & 1          & 1                      & 1                     & 0                          & 0                                      & $-1$                                  \\ \hline
    $11n_{19}$    & 1          & 1                      & 1                     & 0                          & 0                                      & $-1$                                  \\ \hline
    $11n_{20}$    & 0          & 0                      & 0                     & 0                          & 0                                      & 0                                     \\ \hline
    $11n_{38}$    & 0          & 0                      & $-1$                  & 1                          & 1                                      & 1                                     \\ \hline
    $11n_{57}$    & 1          & 2                      & 1                     & 0                          & 0                                      & $-1$                                  \\ \hline
    $11n_{61}$    & 1          & 1                      & 0                     & 0                          & 0                                      & 0                                     \\ \hline
    $11n_{70}$    & 1          & 1                      & 1                     & 0                          & 0                                      & $-1$                                  \\ \hline
    $11n_{79}$    & 0          & 0                      & 0                     & 0                          & 0                                      & 0                                     \\ \hline
    $11n_{96}$    & 0          & 0                      & $-2$                  & 2                          & 2                                      & 2                                     \\ \hline
    $11n_{102}$    & 0          & 0                      & $-1$                 & 1                          & 1                                      & 1                                     \\ \hline
     $11n_{104}$    & 1                          & 1                                      & 1  & 0          & 0                      & $-1$                                                  \\ \hline       
    $11n_{111}$    & 0          & 0                      & 0                    & 1                          & 1                                      & 0                                     \\ \hline
    $11n_{135}$    & 0          & 0                      & $-1$                 & 1                          & 1                                      & 1                                     \\ \hline
    \end{tabular}
    \caption{The involutive concordance invariants $\underline{V}_0$ and $\overline{V}_0$ for the 10- and 11-crossing knots of interest are shown alongside $V_0$ for comparison.}
    \label{tab:invol}
\end{table}

\FloatBarrier

\subsection{Amendments to the literature}\label{subsec:amend}

In the process of studying the Heegaard diagrams for some of the 11-crossing knots of interest, we discovered a few errors in \cite[Section 3.7.3]{Racz}. In Section 3.7.3, Racz lists parameterizations for the Heegaard diagrams of non-Montesino (1,1)-knots up to 12 crossings. The entries for the knots $12n_{404}$, $12n_{579}$, and $12n_{749}$ do not produce valid diagrams. We found a correct parameterization for $12n_{404}$ by elimination the possibilities based on the Alexander polynomial of the knot. The correct result for $12n_{749}$ is given in \cite[Section 6.2] {rasmussen}. \color{black} Table \ref{tab:params} summarizes these two corrections.

\begin{table}[htp]
    \centering
    \begin{tabular}{|c|c|}
    \hline
    \textbf{Knot} & $\bm{(k,r,c,s)}$ \\ \hline
    $12n_{404}$   & $(14,7,-7,1)$ \\ \hline
    $12n_{749}$   & $(7,3,-3,4)$ \\ \hline
    \end{tabular}
    \caption{Valid parameterizations for the Heegaard diagrams of $12n_{404}$ and $12n{749}$}
    \label{tab:params}
\end{table}


\bibliographystyle{amsalpha}
\bibliography{bib}

\end{document}